\documentclass[review]{elsarticle}

\usepackage{tikz}
\usepackage{pgfplots}

\pgfmathsetmacro{\width} {7cm} 
\pgfmathsetmacro{\height} {7cm}

\usepackage{theorem}
\usepackage{amssymb,amsmath}

\usepackage{mathptmx}       
\usepackage{helvet}         
\usepackage{courier}        
\usepackage{type1cm}        
\usepackage{makeidx}         
\usepackage{graphicx}        
\usepackage{multicol}        
\usepackage[bottom]{footmisc}
\usepackage{subfigure}
\usepackage{booktabs}

\usepackage{bbm}
\usepackage{bm}

\usepackage{epsfig} 
\usepackage{epstopdf}
\usepackage{color}
\usepackage{colordvi}
\usepackage[ntheorem]{empheq} 

\newcommand{\innerp}[3]         { ( {#1},{#2} )_{#3} }                      
\newcommand{\norm}[2]         { \| {#1} \|_{#2} }                      
\newcommand{\bfm}[1]             { \mathbf{#1}     }             %
\newcommand{\qq}             { \bfm{q}     }             
\newcommand{\rr}             { \bfm{r}     }             
\newcommand{\Nabla}       { \boldsymbol{\nabla} }   
\newcommand{\SCZO}            { C^0(\Omega) }                    
\newcommand{\SLTO}            { L^2(\Omega) }                       
\newcommand{\SLTOvec}            { [L^2(\Omega)]^2 }                   
\newcommand{\SHOP}            { H^1(\Ph) }                              
\newcommand{\SHdivP}            { H(\text{div},\Ph) }            
\newcommand{\VVK}            { {V(K_m)} }              
\newcommand{\VV}            { {V(\Ph)} }        
\newcommand{\UUU}            { U(\Omega) }           
\newcommand{\UUUT}            { U(\Omega_T) }           
\newcommand{\UUUh}            { U^h(\Omega) }           
\newcommand{\UUUTh}            { U^h(\Omega_T) }           
\newcommand{\VVh}            { V^*(\Ph) }          
\newcommand{\SHOK}            { H^1(K_m) }                        
\newcommand{\SLTOT}            { L^2(\Omega_T) }                    %
\newcommand{\SHOOT}            { H^1(\Omega_T) }                   
\newcommand{\SHdivO}            { H(\text{div},\Omega) }             
\newcommand{\SHdivK}            { H(\text{div},K_m) }            
\newcommand{\SHMHGout}            { H^{-1/2}(\Gout) }          
\newcommand{\SHHdK}            { H^{1/2}(\dKm) }         
\newcommand{\SHmHdK}            { H^{-1/2}(\dKm)  }       
\newcommand{\SLOinf}            { L^\infty(\Omega)}          

\newcommand{\DD}            { \epsilon  }              
\newcommand{\bb}             { {\bfm{b}}   }             
\newcommand{\xx}             { \bfm{x}     }             
\newcommand{\vn}             { \bfm{n}     }             

\newcommand{\Gout}        { \Gamma_{out} }          
\newcommand{\Ph}            { \mathcal{P}_h }
\newcommand{\Kep}           { K_m \in \Ph}
\newcommand{\dKm}           { \partial K_m }

\newcommand{\dx}              { \; {\rm d} \bfm{x}   }                
\newcommand{\dss}             { \, {\rm d} s   }                          
\newcommand{\summa}[2]        { \overset{#2}{\underset{#1}{\sum}} } 
\newcommand{\supp}[1]         { \underset{#1}{\sup} \, }        

\newcommand{\wwwh}              { \bfm{w^*} }                      
\newcommand{\vvh}              { v^* }                   
\newcommand{\vvi}              { \tilde{e}^i }                 
\newcommand{\wwwi}              { \mathbf{\tilde{E}^i}}          
\newcommand{\www}             { \bfm{w} }                      
\newcommand{\nn}              { \bfm{n} }                      
\newcommand{\isdef}           { \overset{\text{def}}{=} } 
\newcommand{\ds}              { \displaystyle }   

\newcommand{\mfu}{u} 
\newcommand{\mfv}{v}
\newcommand{\mfU}{U(\Omega)}
\newcommand{\mfV}{V(\Omega)}
\newcommand{\mfUh}{U^h(\Omega)}
\newcommand{\mfVh}{V^h(\Omega)}
\newcommand{\mfb}{\text{b}}

\newcommand{\mfl}{F}

\newcommand{\mfL}{\mathcal{L}}

\DeclareMathOperator*{\argmin}{arg\,min}
\empheqset{box=\bigfbox}

\theoremstyle{plain}
\theoremheaderfont{\bfseries \upshape}
\newtheorem{thm}{Theorem}[section]

\newtheorem{rem}{Remark}[section]

\newtheorem{prp}{Proposition}[section]

\usepackage[margin=1in]{geometry}

\begin{document}

\begin{frontmatter}
 \title{Automatic Variationally Stable Analysis for Finite Element Computations: Transient Convection-Diffusion Problems}
 
 \author[1,4,5]{Eirik Valseth\corref{cor1}}
\ead{eirik.valseth@nmbu.no}
   
 \author[2]{Pouria Behnoudfar}

\author[1]{Clint Dawson}

\author[3]{Albert Romkes}

 \cortext[cor1]{Corresponding author}

\address[1]{Oden Institute for Computational Engineering and Sciences, University of Texas at Austin, Austin, TX 78712, USA}

\address[4]{The Department of Data Science, The Norwegian University of Life Science, Drøbakveien 31, Ås 1433, Norway }

\address[5]{Department of Scientific Computing and Numerical Analysis, Simula Research Laboratory, Kristian Augusts gate 23, Oslo, 0164, Norway }

\address[2]{Mineral Resources, Commonwealth Scientific and Industrial Research Organisation (CSIRO), Kensington, Perth, WA 6152, Australia}

\address[3]{Department of Mechanical Engineering, South Dakota School of Mines \& Technology, Rapid City, SD 57701, USA \\[0.1in] \rm{We dedicate this work to Prof. Lsezek F. Demkowicz on the occasion of his 70th birthday. }}


\biboptions{sort&compress}


%
%
%
\begin{abstract}
We present an application of stable finite element (FE) approximations of convection-diffusion initial boundary value 
problems (IBVPs) using a weighted least squares FE method, the automatic variationally stable finite element (AVS-FE) method~\cite{CaloRomkesValseth2018}. 
The transient convection-diffusion problem leads to issues in classical FE methods as the 
differential operator can be considered a singular perturbation in both space and time. 
The stability property of the AVS-FE method, allows us significant 
flexibility in the construction of FE approximations in both space and time. 
Thus, in this paper, we take two distinct approaches to the FE discretization of the convection-diffusion problem: $i)$ considering a space-time approach in which the temporal discretization 
is established using finite elements, and  $ii)$ a method of lines approach in which we 
employ the AVS-FE method in space whereas the temporal domain is discretized using the
generalized-$\alpha$  method.
We also consider another space-time technique in which the temporal direction 
is partitioned, thereby leading to finite space-time "slices" in an attempt to  reduce the computational cost of the space-time discretizations.

We present numerical verifications for these approaches, including numerical asymptotic convergence 
studies highlighting optimal convergence properties.
Furthermore, in the spirit of the discontinuous Petrov-Galerkin (DPG) method by Demkowicz and Gopalakrishnan
\cite{Demkowicz4, Demkowicz2, Demkowicz3, Demkowicz5, Demkowicz6}, the  AVS-FE method 
also leads to readily available \emph{a posteriori} error estimates through a Riesz representer of the 
residual of the AVS-FE approximations. Hence, the norm of the resulting local restrictions 
of these estimates serve as error indicators in both space and time for which we present multiple numerical verifications in mesh adaptive strategies. 
\end{abstract}

\end{frontmatter}

\section{Introduction}
\label{sec:introduction}

Transient BVPs are commonplace in engineering applications and to date still pose 
significant challenges in numerical analysis and numerical modeling.
Time dependency in many BVPs, such as the heat equation, involve partial derivatives 
of the trial variable with respect to time  and leads to numerical instabilities unless careful considerations are taken. The reason being that the time derivative
is a convective transport term, i.e., transient problems may lead to unstable discretizations, particularly in the FE context. 
Additionally, the target problem of convection-diffusion also result in numerical 
instabilities in its spatial discretizations which lead to the development of the AVS-FE method in~\cite{CaloRomkesValseth2018}. 
To overcome the stability issues in both space and time
we propose two distinct approaches employing the AVS-FE method. First, we take a space-time 
approach in which space and time are discretized directly considering time an 
additional dimension using the AVS-FE method. Second, we consider a 
method of lines to decouple the computations in space and time and employ a generalized
$\alpha$ method for the temporal discretization~\cite{deng2019high,behnoudfar2019higher,chung1993time}. 

The use of space-time FE methods remains attractive as the approximations are
standard FE approximations and therefore inherit attractive features 
of FE methods such as \emph{a priori} and \emph{a posteriori} error estimation and
mesh adaptive strategies. Examples of space-time FE methods can be found in,
e.g.,~\cite{Hughes1996,hughes1988space,aziz1989continuous}. 
The AVS-FE method~\cite{CaloRomkesValseth2018} being  
stable for any differential operator is therefore a prime candidate for 
space-time FE discretizations. 
Its stability property is a consequence of the philosophy of the DPG method in which the test space 
consist of functions that are computed on-the-fly from Riesz representation problems~\cite{Demkowicz4, Demkowicz2, Demkowicz3, Demkowicz5, Demkowicz6}. 
In~\cite{valseth2020Cahn}, the AVS-FE method is successfully employed in space and time for the Cahn-Hilliard BVP. The goal here was the extension of the AVS-FE method to a nonlinear 
BVP as well as an initial verification of AVS-FE space-time solutions.
Similarly, in~\cite{VALSETH2020113297}, the AVS-FE method is employed for space time solutions of a 
nonlinear transient wave propagation problem, the Korteweg de-Vries equation.
Furthermore, its built-in 
\emph{a posteriori} error estimate and their corresponding error indicators can be directly applied to drive adaptivity. 
The DPG method has been successfully applied to
several transient problems, e.g., convection-diffusion and the Navier-Stokes equations~\cite{ellis2014space,ellis2016robust,roberts2015discontinuous}. 
These space-time formulations are available in the DPG FE code Camellia of Nathan Roberts~\cite{roberts2014camellia}. Recent efforts in DPG methods for transient problems  include the use of optimal testing in time, see e.g.,~\cite{munoz2021dpg,munoz2022error}

Alternatively, the method of lines can be employed to decouple the discretization of space and time where the spatial dimension is discretized  to obtain a semi-discrete system. Then, using a time integrator, the discretization of the temporal domain subsequently results in a fully discrete system of equations.  Here, we employ the AVS-FE method in space and 
the generalized-$\alpha$ method in time. Chung and Hulbert introduced the generalized-$\alpha$ method in~\cite{chung1993time} to solve hyperbolic problems and extended it to parabolic differential equations such as Navier-Stokes equations in~\cite{jansen2000generalized}. The method provides second-order accuracy in the temporal domain as well as unconditional stability. Although the method allows us to  control the numerical dissipation in the high-frequency regions, it delivers adequately accurate results in low-frequency domains. Introduction of a user-defined parameter provides this control and includes the HHT-$\alpha$ method of Hilber, Hughes, Taylor~\cite{hilber1977improved} and the WBZ-$\alpha$ method of Wood, Bossak, and Zienkiewicz~\cite{wood1980alpha}.

In the following, we introduce the AVS-FE method for transient BVPs by taking the two distinct 
approaches introduced above. In Section~\ref{sec:avs-fe} we introduce our model problem 
and notations in addition to a review of the AVS-FE methodology and present 
the AVS-FE weak formulation to be used. In this section we also present the 
discretization of the weak form, an alternative
saddle point structure of the AVS-FE method, and its built-in \emph{a posteriori} error estimate. 
In Section~\ref{sec:time_disc} we present the time discretization techniques: the method of lines using AVS-FE method in space and generalized-$\alpha$ method in time is presented in Section~\ref{sec:method_lines}; and the space-time AVS-FE method in Section~\ref{sec:space_time_elems}. 
Results from numerical verifications for numerous PDEs and applications are presented in Section~\ref{sec:numerical_results}. Finally, we draw conclusions and discuss potential directions of future work 
in Section~\ref{sec:conclusions}.

\section{The AVS-FE Methodology}  
\label{sec:avs-fe}
The AVS-FE method~\cite{CaloRomkesValseth2018} allows us to compute  
stable FE approximations to BVPs for any differential operator, provided its kernel is trivial and the computations of optimal test functions is sufficiently accurate. 
In this section we introduce our model problem and briefly review the AVS-FE method, 
a thorough introduction can be found in~\cite{CaloRomkesValseth2018}.

\subsection{Model Problem and Notation}
\label{sec:model_and_notation}
Let $\Omega \subset \mathbb{R}^N$, $N \leq 2$ be an open bounded domain with Lipschitz boundary $\partial \Omega$ and outward unit normal vector $\vn$, and let $T$ be the final time. Then, define $\Omega_T = \Omega \times (0,T)$ to be the space time domain  which is open and bounded with a  Lipschitz
boundary $\partial \Omega_T = \overline{ \Gamma_{in} \cup \Gamma_{out} \cup \Gamma_{0} \cup \Gamma_{T} }$.
$\Gamma_{in}$ and $\Gamma_{out}$ are the in and outflow boundaries, respectively, and $\Gamma_{0}$ and $\Gamma_{T}$ are the initial and final time boundaries, respectively.
The transient model problem is therefore the following 
linear convection-diffusion IBVP:
\begin{equation} \label{eq:transient_conv_diff_BVP}
\boxed{
\begin{array}{l}
\text{Find }  u  \text{ such that:}    
\\[0.05in] 
\qquad 
\begin{array}{rcl}
\ds \frac{\partial u }{\partial t} -\Nabla \cdot(\DD \Nabla u) \, + \,
 \bb \cdot \Nabla u  & = & f, \quad \text{ in } \, \Omega_T, 
 \\[0.05in]
 \qquad u &  = & u_{in}, \quad \text{ on } \, \Gamma_{in} , 
 \\
 \qquad \DD\Nabla u \cdot \vn & = & g, \quad \text{ on } \,  \Gamma_{out},
  \\
 \qquad u &  = & u_{0}, \quad \text{ on } \, \Gamma_{0} , 
 \end{array}
 \end{array}
}
\end{equation}
where $\DD\in\SLOinf$ denotes the isotropic diffusion parameter; $\bb\in\SLTOvec$ the convection coefficient; $f\in\SLTO$ the source function; and $g\in\SHMHGout$ the Neumann boundary data.
Note that the gradient operator $\Nabla$ refers to the spatial gradient operator, e.g., $\Nabla (\cdot) = \{\frac{\partial (\cdot) }{\partial x},\frac{\partial (\cdot) }{\partial y}  \}^{\text{T}}$.

\subsection{Weak Formulation and FE Discretization}
\label{sec:weak_form}

We omit the full derivation of the weak formulation here and mention key points only.
The derivation of a weak formulation for the AVS-FE method is shown in, e.g.~\cite{CaloRomkesValseth2018}.
To establish a weak formulation of~\eqref{eq:transient_conv_diff_BVP}, we need a regular partition 
$\Ph$ of $\Omega_T$ into elements $K_m$, such that:
\begin{equation}
\notag
\label{eq:domain}
  \Omega_T = \text{int} ( \bigcup_{\Kep} \overline{K_m} ).
\end{equation}
%
%
We introduce a flux variable  
$\qq =\DD\Nabla u$, and recast~\eqref{eq:transient_conv_diff_BVP} as a system 
of (distributional) first-order PDEs:
\begin{equation} \label{eq:conv_diff_BVP_first_order}
\boxed{
\begin{array}{l}
\text{Find }  (u,\qq) \in \SHOOT\times(\SHdivO \times L^2(0,T) ) \text{ such that:}    
\\[0.05in] 
\qquad 
\begin{array}{rcl}
\ds   \Nabla u - \frac{1}{\DD} \qq  & =  & 0, \quad \text{ in } \, \Omega_T, 
  \\
\ds \frac{\partial u }{\partial t}   -\Nabla \cdot \qq \, + \,
 \bb \cdot \Nabla u  & = & f, \quad \text{ in } \, \Omega_T, 
 \\[0.025in]
 \qquad u &  = & u_{in}, \quad \text{ on } \, \Gamma_{in} , 
 \\
 \qquad \qq \cdot \vn & = & g, \quad \text{ on } \,  \Gout,
  \\
 \qquad u &  = & u_{0}, \quad \text{ on } \, \Gamma_{0}.
 \end{array}
 \end{array}
}
\end{equation}
Note that the flux variable $\qq$ depends on time but only has the same number of components as 
the dimension of $\Omega$ and in the weak enforcement of the PDE, it 
belongs to $\SHdivO \times L^2(0,T)$. To make notation more compact 
in the following, we write $\SHdivO$ for $\SHdivO \times L^2(0,T)$, and analogously for its broken 
counterpart.

To derive the AVS-FE weak formulation, we enforce the PDEs~\eqref{eq:conv_diff_BVP_first_order} 
weakly on each element $\Kep$, apply integration by parts to
shift all derivatives to the test functions except the time derivative. After subsequent summation of the local contributions we arrive at the global variational formulation:  
\begin{equation} \label{eq:weak_form}
\boxed{
\begin{array}{ll}
\text{Find } (u,\qq) \in \UUUT & \hspace{-0.05in} \text{ such that:}
\\[0.05in]
 &  \quad B((u,\qq),(v,\www)) = F((v,\www)), \quad \forall (v,\www)\in \VV, 
 \end{array}}
\end{equation}
In~\eqref{eq:weak_form}, the bilinear form, $B:\UUUT\times\VV\longrightarrow \mathbb{R}$,
and linear functional, $F:\VV\longrightarrow \mathbb{R}$, are defined:
\begin{equation} \label{eq:B_and_F}
\begin{array}{c}
B((u,\qq),(v,\www)) \isdef
\ds \summa{\Kep}{}\biggl\{ \int_{K_m}\bigl[  \,  -  u\, \Nabla \cdot \www_m \, - \frac{1}{\DD} \qq \cdot \www_m \, 
+ \, \frac{\partial u }{\partial t} v_m \, +  \, \qq \cdot \Nabla v_m \, - \,
 (\bb \cdot \Nabla v_m) \, u  \bigr] \dx\biggr.
 \\[0.1in]
  \ds 
 + \oint_{\dKm} \biggl[ (\bb \cdot \vn) \, \gamma^m_0(u)  \gamma^m_0(v_m) +   \gamma^m_\nn(\www_m) \, \gamma^m_0(u)   -\gamma^m_\nn(\qq) \, \gamma^m_0(v_m) \, \biggr] \dss  
 \biggr\} ,
 \\[0.15in]
 F((v,\www)) \isdef   \ds \summa{\Kep}{}\int_{K_m} f\,v_m \dx,
 \end{array}
\end{equation}
where the \emph{continuous} trial and \emph{broken} test function spaces, $\UUUT$ and $\VV$,
 are defined as follows:
\begin{equation}
\label{eq:function_spaces}
\begin{array}{c}
\UUUT \isdef \biggl\{ (u,\qq)\in \SHOOT\times\SHdivO: \; u|_{\Gamma_{0}} = u_0, \, u|_{\Gamma_{in}} = u_{in}, \, \qq \cdot \vn|_{\Gamma_{out}}  =  g \biggr\},
\\[0.15in]
\VV \isdef \biggl\{ (v,\www)\in \SHOP\times\SHdivP \biggr\}.
\end{array}
\end{equation}
The broken Hilbert spaces are defined:
\begin{equation}
\label{eq:function_spaceH1}
\begin{array}{c}
\SHOP \isdef \biggl\{ v \in \SLTOT: \; v_m \in \SHOK, \; \forall \Kep   \biggr\},\\
\SHdivP \isdef \biggl\{ v \in \SLTOvec: \; \www_m \in \SHdivK, \; \forall \Kep   \biggr\},
\end{array}
\end{equation}
and the norms on these spaces $\norm{\cdot}{\UUUT}:  \UUUT \!\! \longrightarrow \!\! [0,\infty)$ and $\norm{\cdot}{\VV}: \VV\! \! \longrightarrow\! \! [0,\infty)$ are defined as follows:
\begin{equation}
\label{eq:broken_norms}
\begin{array}{l}
\ds \norm{(u,\qq)}{\UUUT} \isdef \sqrt{\int_{\Omega} \biggl[ \Nabla u \cdot \Nabla u + u^2   + (\Nabla \cdot \qq)^2+\qq \cdot \qq\biggr] \dx }.
\\[0.2in]
\ds   \norm{(v,\www)}{\VV} \isdef \sqrt{\summa{\Kep}{}\int_{K_m} \biggl[  h_m^2 \Nabla v_m \cdot \Nabla v_m + v_m^2   + h_m^2  (\Nabla \cdot \www_m)^2 + \www_m \cdot \www_m\biggr] \dx }.
 \end{array}
\end{equation}
The operators $\gamma^m_0: \SHOK: \longrightarrow \SHHdK$ and $\gamma^m_\nn:\SHdivK \longrightarrow \SHmHdK$ denote the trace and normal trace operators on $K_m$.

\begin{rem} \label{rem:weighted_norm}
With the regularities of trial and test spaces in the weak form (see~\eqref{eq:function_spaces} and~\eqref{eq:function_spaceH1}), it is always possible to integrate back to the trivial weak form in which the PDEs~\eqref{eq:conv_diff_BVP_first_order} are enforced weakly. Hence, the weak form~\eqref{eq:weak_form} represents a first-order system least squares (FOSLS)~\cite{bochevLeastSquares} 
weak form in which the derivatives are shifted to the test functions.  However, our choice of norm on the test 
space $\norm{\cdot}{\VV}$~\eqref{eq:broken_norms} is stronger than the $L^2$ norm used 
in FOSLS. 
This norm is in fact equivalent to the standard $L^2$ norm: 
\begin{equation} \label{eq:weaker_stronger_norm}
 C_1 \norm{(v,\www)}{\VV} \leq \norm{(v,\www)}{\SLTOT} \leq C_2 \norm{(v,\www)}{\VV},
\end{equation}
where the upper bound is obvious and the lower bound holds on quasi-uniform meshes due to an inverse inequality. Hence, $C_1$ is mesh dependent.

Our reasoning for this choice is based on numerical evidence suggesting robust and stable behavior  
for limiting cases of convection dominated diffusion and near incompressibility in linear elasticity (see, e.g.,~\cite{CaloRomkesValseth2018}). 
Stronger norms is also shown to have positive effects in the minimum residual technique introduced in~\cite{calo2019adaptive}. We also point out that the choice makes sense by considering that the order of magnitude of the terms in the norm definition~\eqref{eq:broken_norms} are identical due to the scaling by the element diameter. Furthermore, from an engineering point-of-view, had the test functions had units, all terms would be of the same unit 
due to this scaling.
\end{rem}

The bilinear form and linear functional in~\eqref{eq:B_and_F} differs from the ones 
presented in~\cite{CaloRomkesValseth2018} due to the 
term $\frac{\partial u}{\partial t}$ and the application of integration by parts to all terms involving spatial derivatives. 
This weak formulation~\eqref{eq:weak_form} represents a DPG formulation as the 
test space is broken and continuity of the trial space is a results of the definition of its subspaces. 
setting. 
In the following we review important points of the AVS-FE method and for the sake of simplicity, consider the case with homogeneous Dirichlet boundary conditions  ($u|_{\partial \Omega_T} = 0$) which are enforced strongly in the trial space 
$\UUUT$.

A key point in the well posedness of least-squares, AVS-FE and DPG weak formulations and FE discretizations is the existence and use of 
an equivalent norm on the trial space $\UUUT$.
Since the kernel of $B(\cdot,\cdot)$ is trivial, we introduce the following  \emph{energy norm} $\norm{\cdot}{\text{B}}: \UUUT\longrightarrow [0,\infty)$:
\begin{equation}
\label{eq:energy_norm}
\norm{(u,\qq)}{\text{B}} \isdef \supp{(v,\www)\in \VV\setminus \{(0,\mathbf{0})\}} 
     \frac{|B((u,\qq),(v,\www))|}{\norm{(v,\www)}{\VV}}.
\end{equation}
As in the DPG method,  the energy norm of $(u,\qq) \in \UUUT$ can be identified by functions $(p, \rr) \in \VV$ that are solutions of the following Riesz representation problem:
\begin{equation} \label{eq:riesz_problem}
\begin{array}{rcll}
\ds \left(\, (p,\rr),(v,\www) \, \right)_\VV &  \! \! =  \! & B(\,(u,\qq),(v,\www) \, ),& \, \forall (v,\www)\in\VV, 
\end{array}
\end{equation}
where $\left(\, \cdot, \cdot \, \right)_\VV:\; \VV\times\VV \longrightarrow \mathbb{R}$, is the inner product on $\VV$ defining the norm $\norm{\cdot}{\VV}$ in~\eqref{eq:broken_norms}.
Due to the Riesz representation problem~\eqref{eq:riesz_problem} we can establish the 
equivalence between the energy norm of trial functions $(u,\qq)\in\UUUT$, and the norm of the Riesz representers $(p,\rr) \in \VV$:
\begin{equation}
\label{eq:norm_equivalence}
\norm{(u,\qq)}{\text{B}} = \norm{(p,\rr)}{\VV}.
\end{equation} 
The well posedness in terms of the energy norm is essentially an assumption of this method as evident by its definition~\eqref{eq:energy_norm}. Note that in the FOSLS, the test space is $L^2$ and the solution of the Riesz problem is trivial: 
\begin{equation}
\label{eq:fosls_riesz}
p = \frac{\partial u }{\partial t}   -\Nabla \cdot \qq \, + \,
 \bb \cdot \Nabla u, \; \; \rr = \Nabla u - \frac{1}{\DD} \qq.
\end{equation} 
Hence, the energy norm can be exactly identified by the definition of the weak form and the analyses presented in~\cite{bochevLeastSquares} can be applied.

Next, we present a brief review of the AVS-FE spatial discretization, the 
discretization of the time domain is presented separately in Section~\ref{sec:time_disc}.
Hence, we suppress notation related to time dependency in this and the following section. 
To establish FE approximations $(u^h,\qq^h)$ of $(u,\qq)$ the AVS-FE method follows the 
classical FE method and represents
the FE approximations $u^h$ and $\qq^h$  as linear combinations of 
 basis functions and their corresponding degree of freedom.
Proper choices of bases are, e.g., continuous polynomials for the base variable $u^h \in P^p(\Omega)$ and Raviart-Thomas polynomials for the flux $\qq^h \in RT^p(\Omega)$.

As the test space  $\VV$ is broken, the test functions are to be piecewise 
discontinuous and are constructed by employing the DPG philosophy~\cite{Demkowicz4, Demkowicz2, Demkowicz3, Demkowicz5, Demkowicz6,niemi2013automatically}.
Hence, each basis function in the trial space $\UUUh$ is paired with a (vector valued)
test function.  In the same spirit as $(p,\rr)$ are the Riesz representers of $(u,\qq)$ in~\eqref{eq:riesz_problem}, $(\vvi,\wwwi)$ are the Riesz representers of the basis functions $(e^i, (E_x^j(\xx), E^k_y(\xx)))$ through~\eqref{eq:riesz_problem}.
Clearly,~\eqref{eq:riesz_problem} is of infinite dimension and must be approximated.
Due to the broken nature of $\VV$ we can solve these local problems in a decoupled fashion element-by-element by computing 
piecewise polynomial approximations to the optimal test functions. 

Finally, the discretization governing the FE approximation $(u^h,\qq^h)\in\UUUh$ is:
\begin{equation} \label{eq:discrete_form}
\boxed{
\begin{array}{ll}
\text{Find } &  (u^h,\qq^h) \in \UUUh \; \text{ such that:}
\\[0.1in]
 &   B((u^h,\qq^h),(\vvh,\wwwh)) = F((\vvh,\wwwh)), \quad \forall (\vvh,\wwwh)\in \VVh, 
 \end{array}}
\end{equation}
where the finite dimensional subspace of test functions $\VVh\subset\VV$ is spanned by the numerical 
approximations of the test functions.
This discrete problem is guaranteed to admit stable FE approximations if the test functions are computed 
with sufficient accuracy, which in the is equivalent to the existence (local) DPG Fortin 
operators~\cite{nagaraj2017construction,demkowicz2020construction}. 
In the current work, we do not perform the construction of the Fortin operators, but  perform  
a numerical study in Section~\ref{sec:optimal_res} to verify stability for the canonical AVS-FE 
choice of using discontinuous polynomials of the same degree as the continuous trial space.

\subsection{Saddle Point Problem}
\label{sec:saddle_point_problem}
The AVS-FE discretization~\eqref{eq:discrete_form} can be implemented in existing FE
software by redefining routines that compute the element stiffness matrices. However, 
in several commonly used FE solvers, such as FEniCS~\cite{alnaes2015fenics} or Firedrake~\cite{rathgeber2017firedrake}, manipulations of the element assembly
routines may not as easily be performed.
Thus, to enable straightforward implementation into these FE solvers, we will introduce an equivalent
interpretation of the AVS-FE method as a global saddle point problem.
We omit several details here and highlight only key features of this interpretation, interested readers are referred to~\cite{demkowicz2014overview} for a complete presentation.

The AVS-FE method is a weighted least squares, or minimum 
residual method,  in the sense that its solution realizes the minimum of a functional 
according to the following principle:
\begin{equation} \label{eq:min_DPG_residual}
\boxed{
\begin{array}{ll}
 \ds u^h = \argmin_{v^h\in \UUUh} \frac{1}{2} \norm{\mathbbmss{B}v^h- \mathbbmss{F}}{\VV'}^2, 
 \end{array}}
\end{equation}
where $\mathbbmss{B}$ and $\mathbbmss{F}$ are operators induced by the bilinear and 
linear forms, respectively.
Due to the Riesz representation problem~\eqref{eq:riesz_problem} and energy norm, we  can relate the norm on the dual space $\VV'$ $\norm{\cdot}{\VV'}$ 
to the energy norm $\norm{\cdot}{\text{B}}$. 
Thus, we can consider a Riesz representer of the approximation error $(u-u^h,\qq-\qq^h)$, which we refer to as an \emph{error representation function}~\cite{Demkowicz2}. This error representation function
$( \hat{e},\hat{\bfm{E}} )$ is then defined as the solution 
of the following weak problem:
\begin{equation} \label{eq:error_rep_problem}
\boxed{
\begin{array}{ll}
\text{Find } (\hat{e},\hat{\bfm{E}} ) \in \VV  \quad \text{such that:}
\\[0.05in]
 \innerp{(\hat{e} ,\hat{\bfm{E}} )}{(v,\www)}{\VV} = \underbrace{ F(v,\www) -  B(\, (u^h,\qq^h),(v,\www)\, )}_{\text{Residual}} \quad \forall \, (v,\www)  \in \VV. 
 \end{array} }
\end{equation}
%
%
The energy norm of $(u-u^h,\qq-\qq^h)$ can be identified by the $\VV$ norm of the error representation function:
\begin{prp}
\label{prp:err_rep_equ}
Let $(u,\qq) \in \UUU$  be the solution of the AVS-FE weak form~\eqref{eq:weak_form} and $(u^h,\qq^h) \in \UUUh$ its corresponding AVS-FE approximation through~\eqref{eq:discrete_form}.
Then, the  the energy norm of 
$(u-u^h,\qq-\qq^h)$ is identical to the $\VV$ norm of $ (\hat{e} ,\hat{\bfm{E}})$:
\begin{equation} \label{eq:error_rep_norm}
\norm{(u-u^h, \qq - \qq^h)}{\rm{B}} = \norm{(\hat{e} ,\hat{\bfm{E}} )}{\VV}.
\end{equation}
\end{prp} 
\emph{Proof}: This proof is known from existing DPG literature (see Section 1 and equation (1.17) in~\cite{Demkowicz6}). The identity is a consequence of the norm equivalence in~\eqref{eq:norm_equivalence}, the definition of the energy norm~\eqref{eq:energy_norm} and the weak problem governing the error representation 
function~\eqref{eq:error_rep_problem}.
\newline \noindent ~\qed

\noindent The norm of approximate error representation function $(\hat{e}_h ,\hat{\bfm{E}}_h )$ is therefore an \emph{a posteriori} error estimate, i.e,
\begin{equation} \label{eq:error_est}
\norm{(u-u^h, \qq - \qq^h)}{\rm{B}} \approx  \norm{(\hat{e}_h ,\hat{\bfm{E}}_h )}{\VV}.
\end{equation}
Furthermore, its local restriction can be computed element-wise as the space $\VV$ is broken to yield
the error indicator:
\begin{equation} \label{eq:error_rep_norm_local}
\eta = \norm{(\hat{e}_h ,\hat{\bfm{E}}_h )}{\VVK}.
\end{equation}
This type of error indicator has been applied with great success to multiple problems 
(see, e.g., \cite{Demkowicz2,Demkowicz6,fuentes2017coupled,calo2019adaptive}), 
and we show several numerical experiments using this indicator for the AVS-FE method in Section~\ref{sec:numerical_results}. It should be noted that this error estimate and the error indicator 
are known to be robust (i.e., bounded above and below) under the assumption of the existence of  
DPG Fortin operators and localizable 
norms~\cite{Demkowicz2,nagaraj2017construction,demkowicz2020construction}.

The minimum residual interpretation allows us to establish the following AVS-FE 
saddle point formulation to which we seek the approximate solution $(u^h,\qq^h)$ under the 
constraint of the error representation function minimizes the residual of the AVS-FE method, see~\eqref{eq:error_rep_problem}:
\begin{equation} \label{eq:saddle_point_problem}
\boxed{
\begin{array}{ll}
\text{Find } (u^h,\qq^h) \in  \UUUh, (\hat{e}_h ,\hat{\bfm{E}}_h ) \in V^h(\Ph)  & \hspace{-0.15in} \text{ such that:}
\\[0.05in]
   \quad \left(\, (\hat{e}_h ,\hat{\bfm{E}}_h),(v^h,\www^h) \, \right)_\VV + B((u^h,\qq^h),(v^h,\www^h)) & =   F(v^h,\www^h), \quad \forall (v^h,\www^h)\in V^h(\Ph),  \\
  \quad B((p^h,\rr^h),(\hat{e}_h ,\hat{\bfm{E}}_h ))& =  0, \quad \forall (p^h,\rr^h) \in \UUUh.
 \end{array}}
\end{equation}
Solution of~\eqref{eq:saddle_point_problem} gives both the AVS-FE solution 
for $(u^h, \qq^h)$ and its error representation functions 
$(\hat{e}_h ,\hat{\bfm{E}}_h )$ in a single global solution step.
This is very convenient as we now have a built-in \emph{a posteriori} error estimate and error indicators immediately upon solving~\eqref{eq:saddle_point_problem}. However, the computational cost of doing so has been shifted from local computations for optimal test functions to the global cost of a larger system of equations. Fortunately, the global nature of~\eqref{eq:saddle_point_problem} allows for very simple 
implementation of the AVS-FE method in readily available FE solvers like FEniCS~\cite{alnaes2015fenics} and Firedrake~\cite{rathgeber2017firedrake}. Note that dropping the weighted derivative terms from the inner product corresponding to the norm $\norm{\cdot}{\VV}$ reduces~\eqref{eq:saddle_point_problem} to a DPG implementation of the first-order system least squares method. Note that the analysis of~\eqref{eq:saddle_point_problem} can be performed using the famous Brezzi theory~\cite{BrezziMixed,brezzi1974existence}. Since the inner product is a coercive linear operator, and the 
bilinear form satisfies a discrete \emph{inf-sup} condition, the saddle point system is also well posed.

\section{Time Discretization}
\label{sec:time_disc}
In the weak formulation~\eqref{eq:weak_form} we have made no assumptions on the type 
of discretization of the time domain. Here, we consider two distinctive cases of time discretization techniques. In both cases the spatial discretizations 
are performed with finite elements and the AVS-FE methodology. 
First, we consider a discretization of the time domain by employing the method of lines 
to decouple the spatial and time discretization and subsequently employing the generalized-$\alpha$ method. Second, the discrete stability property of the 
AVS-FE method allows us to discretize the time domain with finite elements in a space-time approach.

\subsection{Method of Lines}
\label{sec:method_lines}
In this section, we first discuss the method in an abstract setting before proceeding to the 
particular case of the AVS-FE method and generalized-$\alpha$ methods.
To this end, we define two Hilbert spaces $\mfU$ and $\mfV$, and introduce a well-posed weak formulation for a transient BVP, e.g., the convection-diffusion problem of Section~\ref{sec:model_and_notation}:
\begin{equation}
\boxed{
\begin{array}{ll}
\text{Find } \mfu \in \mfU & \hspace{-0.05in}  \text{ such that:} \\[0.05in]
\ds & \mfb(\mfu,\mfv)=\mfl(\mfv), \quad \forall \, \mfv\in \mfV,
 \end{array}}
\label{eq:gen_weak}
\end{equation}
where the bilinear form $\mfb$ contains all spatial and temporal terms.
We then denote by $\mfL$ the time derivative operator, and modify the bilinear form to contain only 
spatial terms and denote it by  $\mfb_h$. To seek approximations of~\eqref{eq:gen_weak} we consider FE polynomial subspaces of $\mfU$ and $\mfV$, i.e., $\mfUh$ and $\mfVh$ and introduce the semi-discrete formulation:
\begin{equation}
\boxed{
\begin{array}{ll}
\text{Find } \mfu^h \in \mfUh & \hspace{-0.05in}  \text{ such that:} \\[0.05in]
\ds & \left(\mfL({\mfu^h}),\mfv^h\right)_{\SLTO}+\mfb_h(\mfu^h,\mfv^h) = \mfl(\mfv^h), \quad \forall \, \mfv^h\in \mfVh,
 \end{array}}
\label{eq:gen_weak_abs}
\end{equation}
where $(\cdot,\cdot)_{\SLTO}$ denotes the $\SLTO$ inner product.
This semi-discrete formulation is assumed to be well-posed.

To advance the solution in time, we consider a uniform partition of the time domain from $t_0 = 0$ to 
the final time $t_N = T$, with  $\tau$ the distance between each step $t_i$. We  compute approximations to $\mfu^h$ at each step using second-order 
accurate generalized-$\alpha$ methods presented in~\cite{chung1993time,jansen2000generalized}.
For  parabolic or first-order hyperbolic problems, the generalized-$\alpha$ method for the transient term $\mfL({\mfu^h})$ in \eqref{eq:gen_weak_abs} is to find  $ \mfu_h^{n+1} \in \mfUh$, such that: 
\begin{equation}\label{eq:par}
(\vartheta_h^{n+\alpha_m} , \, \mfv_h)_{\SLTO}+\mfb_h(\mfu_h^{n+\alpha_f}, \, \mfv_h) =  \mfl^{n+\alpha_f}(\mfv_h),   \quad\forall\, \mfv_h \in \mfV_h,
\end{equation}
where $\mfu_h^n,\, \vartheta_h^n$ are the approximations to $\mfu(., t_n)$ and $\frac{\partial \mfu(.,\,t_n)}{\partial t}$, respectively. The unknowns at time step $n+1$ are updated using the solutions at $n+\alpha_f$ and $n+\alpha_f$ as:
%
%
\begin{equation} \label{eq:mf}
\boxed{
\begin{aligned}
\eta^{n+\alpha_g} & = \eta^n + \alpha_g \delta(\eta^n), \text{ where } \quad \eta = u,\, \vartheta, \,\mfl\quad g=m, \,f,\\
\delta(\eta^n) & = \eta^{n+1} - \eta^n. \\
\end{aligned} }
\end{equation}
Using a Taylor expansion, we have $\mfu^{n+1} = \mfu^n + \tau \vartheta^n + \tau \gamma \delta(\vartheta^n)$ as a linear combination of $\mfu^n, \vartheta^n$ with $\gamma$ guaranteeing second-order accuracy. 
Substitution of the expressions in~\eqref{eq:mf} into~\eqref{eq:par} gives:
\begin{equation}\label{eq:par2}
(\vartheta_h^{n+1} , \, \mfv_h)_{\SLTO}+\mfb_h(\zeta \, \vartheta_h^{n+1}, \, \mfv_h) = (\frac{1}{\alpha_m} l^{n+1},\mfv_h),   \quad\forall\, \mfv_h \in \mfV_h,
\end{equation}
where $\zeta = \frac{\tau \gamma \alpha_f}{\alpha_m}$, and:
\begin{equation}\label{eq:A}
l^{n+1} = \mfl^{n+\alpha_f}+(\alpha_m-1) \, \left(\vartheta_h^n , \, \mfv_h\right)_{\SLTO}+ \tau  \alpha_f(\gamma-1) \, \mfb_h\left(\vartheta_h^n , \, \mfv_h\right)-  \mfb_h\left(\mfu^n_h , \, \mfv_h\right).
\end{equation}
It can be shown that this scheme is formally second order accurate (see~\cite{jansen2000generalized}) if we select:
\begin{equation} \label{eq:gm}
\gamma = \frac{1}{2} + \alpha_m - \alpha_f.
\end{equation}
Finally, to control the numerical dissipation in case of poor spatial resolution, the two parameters $\alpha_m$ and $\alpha_f$ are defined in terms of the spectral radius $\rho_\infty$ corresponding to an infinite time step:
\begin{equation} \label{eq:hfd}
\alpha_m = \frac{1}{2} \Big( \frac{3 - \rho_\infty}{1+\rho_\infty} \Big), \qquad \alpha_f = \frac{1}{1+\rho_\infty}.
\end{equation}

\begin{rem}\label{rem:remark1}
The generalized-$\alpha$ method requires additional initial data for $\vartheta_h^0$. 
This value is obtained by setting  $\alpha_{f}=\alpha_m=n=0$  and solution of~\eqref{eq:par}.
\end{rem}

\begin{rem}\label{rem:remarkro}
	The spectral radius $\rho_\infty$ is a user-defined parameter that provides control on the numerical dissipation such that for $\rho_\infty=1$ there is no dissipation control, and the maximum control is delivered by setting $\rho_\infty=0$. Numerical dissipation can occur for example in the case of poor spatial resolution (for more details, see, \cite{behnoudfar2018variationally}). 
\end{rem}

\subsubsection{Generalized-$\alpha$ and the AVS-FE Method}
\label{sec:conv_diff_alpha}
Having introduced the generalized-$\alpha$ method for a well defined weak formulation, we 
now extend it to the AVS-FE method for our model IBVP of convection-diffusion.
Hence, let us consider the AVS-FE weak formulation~\eqref{eq:weak_form}, and the trial 
and test spaces $\UUU$ and $\VV$ analogous to~\eqref{eq:function_spaces}.
The generalized-$\alpha$ method for the AVS-FE method is:
\begin{equation} \label{eq:par3}
\boxed{
\begin{array}{ll}
\text{Find } (\vartheta_h^{n+1} , \, \qq^{n+1}_h) \in \UUUh   \text{ such that:}
\\[0.05in] 
 \qquad (\vartheta_h^{n+1} , \, \vvh)_{\SLTO}+B_h((\zeta \, \, \vartheta_h^{n+1}, \, \qq_h^{n+1}),(\vvh,\wwwh)) = (\frac{1}{\alpha_m} \ell^{n+1}((\vvh,\wwwh)), \\[0.05in]   \qquad  \qquad  \quad\forall\, (\vvh,\wwwh)\in \VVh, 
 \end{array}}
\end{equation}
where the operators are defined:
\begin{equation} \label{eq:alpha_operators}
\begin{array}{l}
B_h((u,\qq),(v,\www)) \isdef
\ds \ds \summa{\Kep}{}\biggl\{ \int_{K_m}\biggl[  \,  -  u\, \DD\Nabla \cdot \www_m \, - \qq \cdot \www_m \, 
+   \, \qq \cdot \Nabla v_m \, - \,
 (\bb \cdot \Nabla v_m) \, u  \biggr] \dx\biggr.
 \\[0.1in]
  \ds 
 + \oint_{\dKm} \biggl[ (\bb \cdot \vn) \, \gamma^m_0(u)  \gamma^m_0(v_m) +  \gamma^m_\nn(\www_m) \, \gamma^m_0(u)   -\gamma^m_\nn(\qq) \, \gamma^m_0(v_m) \, \biggr] \dss  \biggr\} ,  
\\[0.15in]
 \ell^{n+1}((v,\www)) \isdef 
\ds \summa{\Kep}{} \int_{K_m} (f^{n+\alpha_f} \, v) \dx 
+(\alpha_m-1)\left(\vartheta_h^n , \, v\right)_{\SLTO} \\[0.1in] \ds + \tau  \alpha_f(\gamma-1) \cdot B_h\left( (\vartheta_h^n,\bfm{0}) , \, (v,\www)\right)
-B_h\left((u^n_h,\qq^n_h) , \, (v,\www)\right).
\end{array}
\end{equation}

To establish the solutions to~\eqref{eq:par3} we take the same approach
introduced in Section~\ref{sec:saddle_point_problem} and define a saddle point system similar to~\eqref{eq:saddle_point_problem}. 
The major difference between the "original" weak form~\eqref{eq:weak_form} and the one corresponding 
to the generalized-$\alpha$ method, i.e.,~\eqref{eq:par3} other than the adjusted bilinear and linear forms, is the term $(\vartheta_h^{n+1} , \, \vvh)_{\SLTO}$.
Analogous to the case in Section \ref{sec:saddle_point_problem}, the approximation to~\eqref{eq:par3} is governed by the following minimization problem:
\begin{equation}\label{eq:min_prob}
\boxed{
\begin{array}{l}
\displaystyle  \vartheta_h^{n+1}   = \argmin_{z_h \in \UUUh} \dfrac{1}{2}\| \mathbbmss{l}^{n+1}- \left(\mathbbmss{M}+\zeta \mathbbmss{B}_h\right) \, z_h\|^2_{V_h'},
\end{array}
}
\end{equation}
where the operators $\mathbbmss{B}_h$ and $\mathbbmss{l}^{n+1}$ correspond to the actions of the adjusted forms $B_h$ and $ \ell^{n+1}$, respectively, and $\mathbbmss{M}$ to the new term $(\vartheta_h^{n+1} , \, v)_{\SLTO}$. Thankfully, the Riesz map (induced by the equivalent of 
the Riesz representation problem~\eqref{eq:riesz_problem} for~\eqref{eq:par3})
allows us to relate 
the norm on the dual space $\norm{\cdot}{V_h'}$ to the energy norm on $\UUU$ exactly as in~\eqref{eq:norm_equivalence}. 
Hence, we define the following error representation function: 
\begin{equation} \label{eq:e_h}
\boxed{
\begin{array}{ll}
\text{Find } (\hat{e}^{n+1},\hat{\bfm{E}}^{n+1}) \in \VV  \quad \text{such that:}
\\[0.05in]
\innerp{(\hat{e}^{n+1},\hat{\bfm{E}}^{n+1})}{(v,\www)}{\VV} =  \underbrace{ \ell^{n+1}(v,\www) - (\vartheta_h^{n+1} , \, v)_{\SLTO}+B_h((\zeta \, \, \vartheta_h^{n+1}, \, \qq_h^{n+1}),(v,\www))}_{\text{Residual}}, \\ \hfill \qquad \forall \, (v,\www)  \in \VV. 
 \end{array} }
\end{equation}
which now measures how far we are from the best approximation of $(\vartheta_h^{n+1}, \qq_h^{n+1})$ at the current time step.
In the same fashion as in Section~\ref{sec:saddle_point_problem}, the norm of this function is an \emph{a posteriori} error estimate and its restriction to each $\Kep$ an error indicator.
We finally can introduce the saddle point problem for each time step:
%
%
%
\begin{equation} \label{eq:mix_par}
\hspace{-0.365cm}\boxed{
\begin{array}{ll}
\text{Find } (\vartheta_h^{n+1},\qq_h^{n+1}) \in  \UUUh, (\hat{e}_h^{n+1},\hat{\bfm{E}}_h^{n+1}) \in V_h(\Ph)   \text{ such that:} \\
((\hat{e}_h^{n+1},\hat{\bfm{E}}_h^{n+1}) \, , \, (v_h,\www_h))_{V_h} + ( (\vartheta_h^{n+1},\bfm{0})  , \, (v_h,\www_h))_{\SLTO} +\zeta \cdot B_h( (\vartheta_h^{n+1},\qq_h^{n+1}) , \, (v_h,\www_h)) &=\frac{1}{\alpha_m} \ell^{n+1}((v_h,\www_h)), \\ \hfill  \quad\forall\, (v_h,\www_h) \in V_h(\Ph), \\	((z_h,\rr_h )\, , \,(\hat{e}_h^{n+1},\bfm{0}) )_{\SLTO}+\zeta \cdot B_h((z_h,\rr_h ) \, , \, (\hat{e}_h^{n+1},\hat{\bfm{E}}_h^{n+1}) ) &=  0,\\ \hfill  \quad\forall\, (z_h,\rr_h) \in \UUUh,
 \end{array}}
\end{equation}
where the inner product $( \cdot  , \cdot)_{V_h}$ is defined:
\begin{equation}\label{eq:norm}
((\hat{e}_h^{n+1},\hat{\bfm{E}}_h^{n+1}) \, , \, (v_h,\www_h))_{V_h} =((\zeta \cdot \hat{e}_h^{n+1},\hat{\bfm{E}}_h^{n+1}) \, , \, (v_h,\www_h))_{V(\Ph)}+(\hat{e}_h^{n+1},\bfm{0}) \, , \, (v_h,\www_h))_{\SLTO}.
\end{equation}
Computing $\vartheta_h^{n+1}$ from~\eqref{eq:mix_par}, we  obtain $u_h^{n+1}$ from a Taylor expansion at each time step. The overall procedure requires a  matrix solve at each time step as well as two explicit updates.
Additionally, we maintain the consistency of problem which can be  checked by setting the time step to zero.

Next, we show that our proposed saddle-point problem~\eqref{eq:mix_par} is unconditionally stable in the temporal domain. To achieve this, we must show that our AVS-FE spatial discretization scheme does not alter the unconditional stability of generalized-$\alpha$ method.

\begin{thm}\label{th:1}
	The saddle-point problems in \eqref{eq:mix_par} provides unconditionally stable solutions in temporal domain.
\end{thm}
\emph{Proof:} Our proof relies on established bounds from literature for the generalized-$\alpha$ method \cite{behnoudfar2018variationally}, and reasoning based on the properties of the AVS-FE saddle point problem \eqref{eq:mix_par}.
By applying the generalized-$\alpha$ method on the continuous parabolic problem \eqref{eq:conv_diff_BVP_first_order} to discrete the temporal domain, we obtain a semi-discrete problem that can be written: 
		\begin{equation}
		\begin{bmatrix}
		u^{n+1}\\ \tau \vartheta^{n+1}
		\end{bmatrix}=\Xi
		\begin{bmatrix}
		u^{n}\\ \tau \vartheta^{n}
		\end{bmatrix}+\Pi \, l^{n+\alpha_{f}},
		\end{equation} 
with $\Xi$ and $\Pi $ being a $2\times2$ amplification matrix and a $2\times1$ matrix, respectively. The amplification matrix allows us to write the solution at time step $n+1$ using initial condition and forcing term. The derivation and development of the amplification matrix can be found in, e.g., ~\cite{chung1993time, behnoudfar2019higher, deng2019high}. If the  eigenvalues of this amplification matrix are bounded by one, the method is stable. Hence:
\begin{equation}\label{eq:u}
			\tau \|\vartheta^{n+1}\|_{\SLTO}^2\leq \tau \|\vartheta^{n}\|_{\SLTO}^2+\dfrac{1}{\tau}\|\pi_2\,l^{n+\alpha_{f}}\|^2_{\SLTO},
\end{equation}
where $\pi_2$ is the second component of $\Pi$. Considering the saddle-point problem \eqref{eq:mix_par} with unknown $\vartheta_{h}^{n+1}$, we add the term $\|\vartheta^{n+1}_h-\vartheta^{n+1}\|_{\SLTO}\,$ to the right-side of the inequality \eqref{eq:u} and the inequality still holds. Next, using~$\|\vartheta^{n+1}-\vartheta_{h}^{n+1}\|_{\SLTO}\leq \|\vartheta^{n+1}-\vartheta_{h}^{n+1}\|_{\VV}$, the Cauchy-Schwartz inequality, and the error representation provided by the AVS-FE method, we get: 
		\begin{equation*}
		\begin{aligned}
		\tau \|\vartheta^{n+1}_h\|_{\SLTO}^2&\leq \tau \|\vartheta^{n}\|_{\SLTO}^2+\tau \|\vartheta^{n+1}_h-\vartheta^{n+1}\|_{\SLTO}^2+\frac{1}{\tau}\|\pi_2\, l^{n+\alpha_{f}}\|^2_{\SLTO}\\
		&\leq \tau \|\vartheta^{n}\|_{\SLTO}^2+C\sqrt{\tau} \|(\hat{e}_h^{n+1}\|_{\VV}^2+\frac{1}{\tau}\|\pi_2\, l^{n+\alpha_{f}}\|^2_{\SLTO},\\ \hspace{0.4cm}
		&\leq\tau \|\vartheta^{0}\|_{\SLTO}^2+\sum_{j=0}^{j=N-1} \left(C\sqrt{\tau}\|(\hat{e}_h^{j}\|_{\VV}^2+\frac{1}{\tau}\|\pi_2\, l^{j+\alpha_{f}}\|^2_{\SLTO}\right),\end{aligned}
		\end{equation*}
		where $C>0$ is a constant. 
		Hence, the solution is bounded by the initial solution, forcing, and error representation terms. \newline \noindent ~\qed

\subsubsection{Retrieving initial data}

As pointed out in Remark~\ref{rem:remark1}, we need to retrieve the additional initial data
$\vartheta_h^0$ to solve \eqref{eq:mix_par}. Hence, we set $\alpha_f = \alpha_m  = 0$ and get:
\begin{equation} \label{eq:mix_par0}
\boxed{
\begin{array}{l}
\text{Find } (\vartheta_h^{0},\qq_h^{0}) \in  \UUUh, (\hat{e}_h^{0},\hat{\bfm{E}}_h^{0}) \in V_h(\Ph)   \text{ such that:} \\
\begin{array}{ll}
((\hat{e}_h^{0},\hat{\bfm{E}}_h^{0} \, , \, (v_h,\www_h))_{V_h} + ( (\vartheta_h^{0},\bfm{0})  , \, (v_h,\www_h))_{\SLTO}   &= \ell^0((v_h,\www_h)) -\zeta \cdot b_h( (u^{0},\qq_h^{0}) , \, (v_h,\www_h)) , \\ &  \quad\qquad\forall\, (v_h,\www_h) \in V_h(\Ph), \\
	((z_h,\rr_h )\, , \,(\hat{e}_h^{0},\bfm{0}) )_{\SLTO} &=  0,  \quad\forall\, (z_h,\rr_h) \in \UUUh,
	\end{array}
 \end{array}}
\end{equation}
where $u^0$, $\qq^0$, and $\ell^0((v_h,\www_h))$ correspond to the initial data. To ascertain that the problem for the initial data is well posed~\eqref{eq:mix_par0}, we have the following proposition.
\begin{prp} \label{lem:inf-sup}
Let $(v_h,\www_h) \in V_h$ be arbitrary test functions. Then, $\vartheta_h^{0} \in \UUUh$ exists and 
is unique.
\end{prp}
We omit the proof here as it is trivial to show that $( (\vartheta_h^{0},\bfm{0})  , \, (v_h,\www_h))_{\SLTO}$, i,e, the $\SLTO$ inner product, satisfies the following three properties:
	\begin{itemize}
		\item 	Stability: There exist a constant $C_{\emph{sta}}>0$ independent of the mesh size, such that:
		\begin{equation}\label{eq:infsup_h}
		\inf_{0\neq z_h \in \UUUh} \sup_{0\neq v_h \in V_h} \dfrac{|( (z_h,\bfm{0})  , \, (v_h,\www_h))_{\SLTO}|}{\norm{z_h}{\SLTO}\norm{v_h}{\SLTO}}  \geq C_{\emph{sta}}.
		\end{equation}
		\item Consistency: Employing a similar argument as~\cite{calo2019adaptive} to study the consistency of the saddle-point problem, we can state the consistency as:
		\begin{equation}
		\begin{aligned}
		&( (\vartheta_h^{0},\bfm{0})  , \, (v_h,\www_h))_{\SLTO}= (f^{0},(v_h,\www_h))\\&-\zeta \cdot b_h( (u_h^{0},\qq_h^{0}) , \, (v_h,\www_h)),   \quad\forall\, (v_h,\www_h) \in V_h
		\end{aligned}
		\end{equation} 
		
		\item 	Boundedness: There exists a constant $C_{\emph{bnd}}<\infty$, uniformly with respect to the mesh size, such that:
		\begin{equation}\label{eq:continuity}
		( (z,\bfm{0})  , \, (v_h,\www_h))_{\SLTO} \leq C_{\emph{bnd}} \, \|z\|_{\SLTO} \norm{v_h}{\SLTO}, \quad \forall \, (z, v_h) \in U \times V_{h}.
		\end{equation}
		
	\end{itemize}

See~\cite{di2011mathematical} for details on these conditions.  
\newline \noindent ~\qed

Thus, using \eqref{eq:mix_par0}, we have a stable and adaptive method to find the initial data which is  
critical for the generalized-$\alpha$ method to ensure second-order accuracy in time.

\subsection{Space-Time FE Approach}
\label{sec:space_time_elems}
The use of  FE discretizations for transient problems is commonly avoided due to the inherently
unstable nature of transient problems.
The discretizations must be very 
carefully constructed to achieve discrete stability using the classical FE method. 
However, the stability of the AVS-FE method 
allows us to discretize the entire space-time domain with finite elements in a straightforward manner. Furthermore,  \emph{a posteriori} error estimates and error indicators are immediately available to us as 
error indicators are obtained directly in the saddle point approach of the AVS-FE 
method~\eqref{eq:saddle_point_problem}. 

To establish AVS-FE space-time approximations of weak formulation~\eqref{eq:weak_form} 
or~\eqref{eq:saddle_point_problem},
we pick appropriate discretizations of the space $\SHOOT\times\SHdivO$.  
For $\SHOOT$, the choice is classical FE basis functions that are $C^0$ continuous
 functions in $\Omega_T$ such as Lagrange or Legendre polynomials. For $\SHdivO$, a
conforming choice of basis is, e.g., a Raviart-Thomas basis.
However, as in~\cite{CaloRomkesValseth2018},  we employ  approximations for $\qq^h$ by vector valued $\SCZO$  polynomials for each of its components as this has shown to yield superior 
results for convex domains and sufficiently regular sources. The discretized weak form is therefore: 
\begin{equation} \label{eq:weak_form2}
\boxed{
\begin{array}{ll}
\text{Find } (u,\qq) \in \UUUTh & \hspace{-0.05in} \text{ such that:}
\\[0.05in]
 &  \quad B((u^h,\qq^h),(v^*,\www^*)) = F(v^*,\www^*), \quad \forall v^*,\www^*\in \VVh, 
 \end{array}}
\end{equation}
where the components of $\UUUTh$ are spanned by continuous FE basis functions and $\VVh$ by the optimal test functions.

\subsubsection{Time Slice Approach}
\label{sec:time_slices}
As an alternative to the space-time discretization of the 
full space-time domain $\Omega_T$, in this section we introduce a time slice approach for the AVS-FE method.  
While the space-time approach introduced in the preceding section allows straightforward implementation of the AVS-FE method and its 
"built-in" error indicator can drive mesh adaptive refinements, the large number of degrees of freedom  quickly makes the method intractable. 
In an effort to localize the computational cost of the space-time approach, we propose 
to partition the space-time domain into "space-time slices". The slices can be constructed in a number of ways, from uniformly to a graded mesh structure as considered in~\cite{ellis2016space,ellis2014space} for the DPG method.

\begin{figure}[h!]
\centering
\scalebox{.4}{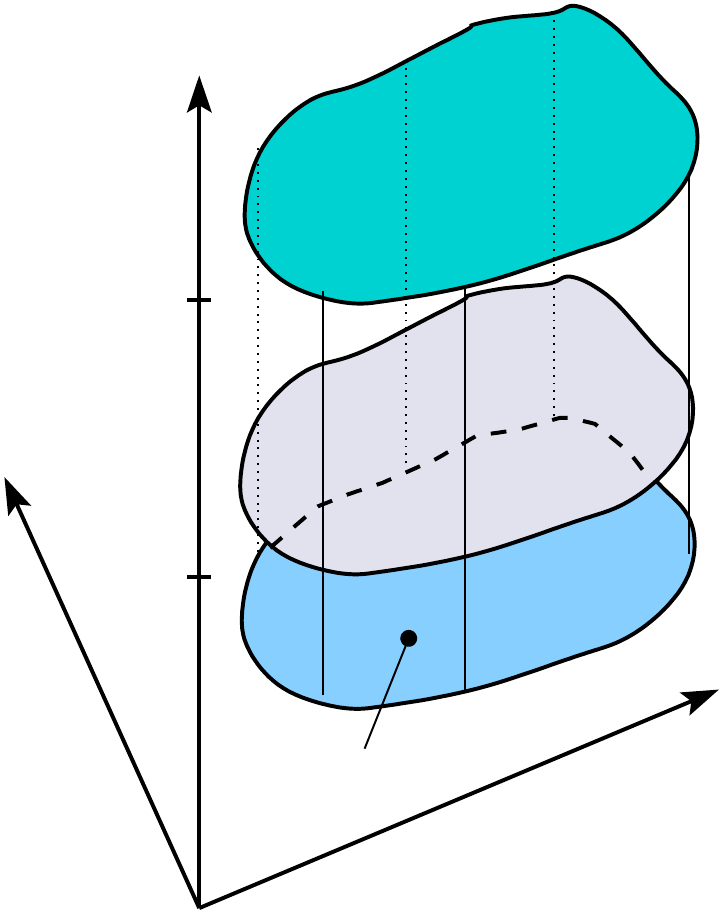 }
\caption{Partition of space-time domain into slices.}
\label{fig:space_time_slices}
\end{figure}

To advance in time, a solution can be 
obtained on a slice which can be transferred to the neighboring slice as an initial 
condition. Hence, we can perform mesh refinements on each slice to ensure the 
complete resolution of any interior or boundary layer (i.e., physical features) before 
proceeding to the next. This is of particular interest in applications in which physical parameters are time dependent leading to widely different solution features as time progresses. 
In Figure~\ref{fig:space_time_slices}, an arbitrary domain $\Omega_T$ is shown and is partitioned  into two space-time slices $\Omega\cap(0,T_{slice})$ and $\Omega\cap(0,T)$. 

\section{Numerical Verifications}
\label{sec:numerical_results}
To conduct numerical verifications, we consider the following form of our model scalar-valued convection diffusion problem~\eqref{eq:transient_conv_diff_BVP} with exact solutiont boundary conditions:
\begin{equation}\label{eq:num_results_homo_conv_diff}  
\begin{array}{rl}
\ds \frac{\partial u }{\partial t} - \epsilon \Delta u +
 \bb \cdot  \Nabla u = f, & \quad \text{ in } \Omega, 
 \\[0.1in]
  u = u_0, & \quad \text{ on } \partial \Omega,  
   \\[0.1in]
  u = u_{intial}, & \quad \text{ on } \partial \Omega \cap \{t=0\},  
 \end{array} 
\end{equation}
where  the coefficient $\epsilon$ is a constant diffusion coefficient.
We first study the effect of approximation degree of the optimal test functions in Section~\ref{sec:optimal_res}.  
Next, we verify the convergence properties of the AVS-FE method for both time discretization schemes in Section~\ref{sec:convergence}. In this section, we also investigate the use of  time slices as well as compare the space-time method to the method of lines with generalized-$\alpha$ time stepping.
Last, in Section~\ref{sec:shock} we  present verifications of a problem with 
both a hyperbolic and a parabolic part, i.e., a transient convection-diffusion
problem. The particular case we investigate corresponds to a  challenging physical application, a shock wave problem.

In all the presented numerical experiments we use the saddle point description in~\eqref{eq:saddle_point_problem} implemented in legacy FEniCS~\cite{alnaes2015fenics} with the latest stable release from Anaconda.  The verifications in which we employ adaptive refinements all use 
the same criterion as in~\cite{calo2019adaptive}, i.e.,
the built-in error indicator~\eqref{eq:error_rep_norm_local}  as well as a D\"orfler marking 
strategy~\cite{dorfler1996convergent} (we pick the parameter $\theta=0.5$) using the approximate energy error computed 
using~\eqref{eq:error_rep_norm}. To solve the system of linear algebraic equations, we use the direct solver MUMPS~\cite{amestoy2006hybrid,amestoy2001fully}.
Also note that in all cases where we report the number of 
 degrees of freedom, we do not include the degrees of freedom for the error representation function 
 in the saddle point systems~\eqref{eq:saddle_point_problem} and~\eqref{eq:mix_par}.
The polynomial degree of approximation used for this error representation function is  
identical to the degree of the trial space with the results in Section~\ref{sec:optimal_res} being the sole exception.

\subsection{Optimal Test Function Resolution}
\label{sec:optimal_res} 

As an initial verification, we perform a study to ensure proper resolution of the optimal test space. 
To this end, we consider the following exact solution: 
\begin{equation}\label{eq:test_res_ex}  
\begin{array}{rl}
\ds u(x,y,t) = e^{-t}\left[x + \frac{e^{ \frac{b_x  }{\epsilon} x}-1}{1-e^{\frac{b_x  }{\epsilon}}}\right]\left[y + \frac{e^{\frac{b_y  }{\epsilon} y}-1}{1-e^{\frac{b_y  }{\epsilon} }}\right], 
 \end{array} 
\end{equation}
from which we establish initial and exact solutiont boundary conditions and a corresponding source term $f$. 
For these studies we consider the moderately convection dominated case with $\epsilon = 0.1$, $\bb = \{1,1 \}$, and select the final time of computation to be $T=0.5s$. 
We consider only the space-time case here and assume that the conclusions apply to the 
generalized-$\alpha$ case as well. Due to the smoothness of the exact solution, we consider continuous
polynomial approximations for both variables of equal order - $p$. The error representation functions 
are then discretized with discontinuous polynomials of order $p+0,1,2,3$, as well as $p-1$ for $p\ge2$. In Table~\ref{tab:test_functs}, these results are presented for linear and quadratic trial functions for two uniform meshes: 6 and 24,576 space-time tetrahedrons, respectively. The results in these table indicate that for linear and quadratic bases for the trial space, the impact of increasing test space degree is vanishing small. We observe the same trend for $p>2$. Note that for $p=2$, we observe satisfactory results for a test space degree $p=1$.
\begin{table}[h]
\centering
\caption{\label{tab:test_functs}  Increasing degree of approximation for the error representation function. }
\begin{tabular}{@{}llll@{}}
\toprule
{$p_{trial}$ \hspace{6mm}} & {$p_{test}$ \hspace{6mm}} & {$\norm{u-u^h}{\SLTOT}$ (coarsest mesh) \hspace{6mm}} & { $\norm{u-u^h}{\SLTOT}$ (finest mesh) \hspace{15mm}}  \\
\midrule \midrule

1 & 1 & 1.1439e-01 & 2.6374e-03    \\
1 & 2 & 1.1439e-01 & 2.6559e-03    \\
1 & 3 & 1.1439e-01 & 2.6576e-03   \\
1 & 4 & 1.1439e-01 & 2.6582e-03   \\
2 & 1 & 6.8837e-02 & 1.4028e-04   \\
2 & 2 & 6.7822e-02 & 1.3770e-04   \\
2 & 3 & 6.8246e-02 & 1.3716e-04   \\
2 & 4 & 6.8277e-02 & 1.3698e-04  \\
2 & 5 & 6.8288e-02 & 1.3695e-04  \\

\bottomrule
\end{tabular}
\end{table}

\subsection{Convergence Studies}
\label{sec:convergence}

To numerically investigate the convergence properties of our methods,
we consider a well-known example of transient convection-diffusion, the Eriksson-Johnson problem \cite{eriksson1993adaptive}. This problem has a known exact solution that satisfies the following form 
of~\eqref{eq:num_results_homo_conv_diff}: 
\begin{equation}\label{eq:conv_diff_eriksson}  
\begin{array}{rl}
\ds \frac{\partial u }{\partial t}  - \epsilon \, \Delta u +
 \frac{\partial u }{\partial x} = f, & \quad \text{ in } \Omega_T.  
 \end{array} 
\end{equation}
Additionally, Dirichlet boundary conditions on $u$, the initial condition on $u$, and the source $f$ are ascertained from the exact solution:
\begin{equation}\label{eq:conv_diff_eriksson_exact}  
\begin{array}{rl}
\ds u_{ex} (\xx) = e^{\ds - l \, t}\, \left( e^{\ds \lambda_1 \, x}\, - e^{\ds \lambda_2 \, x}\,  \right) + \text{cos} (\pi \, y)\, \ds \frac{\ds e^{\ds \delta_2 \, x}\, - e^{\ds \delta_1 \, x} }{\ds e^{\ds - \delta_2 }\, - e^{\ds - \delta_1} }   ,  
 \end{array} 
\end{equation}
where $\ds l = 2$, and:
\begin{equation} \label{eq:eriksson_johnson_params}
\begin{array}{rcl}
\ds \lambda_{1,2} = & \frac{\ds - 1 \pm \sqrt{1 - 4 \, \epsilon \, l}}{\ds -2 \, \epsilon },  \\[0.05in]
\ds \delta_{1,2} = & \frac{\ds  1 \pm  \sqrt{1 + 4 \, \pi^2 \, \epsilon^2}}{\ds 2 \, \epsilon },   \\[0.05in]
 \end{array}
\end{equation}
The problem domain $\Omega_T = ( -1, 0 ) \times ( -0.5,0.5  ) \times ( 0 ,0.5 ) $.
For these studies we consider the moderately convection dominated case of \eqref{eq:conv_diff_eriksson} with $\epsilon = 0.075$.

\begin{figure}[h!]
\subfigure[ \label{fig:asympt_convp1_spacetime} Linear polynomial approximations. ]{\centering
 \includegraphics[width=0.51\textwidth]{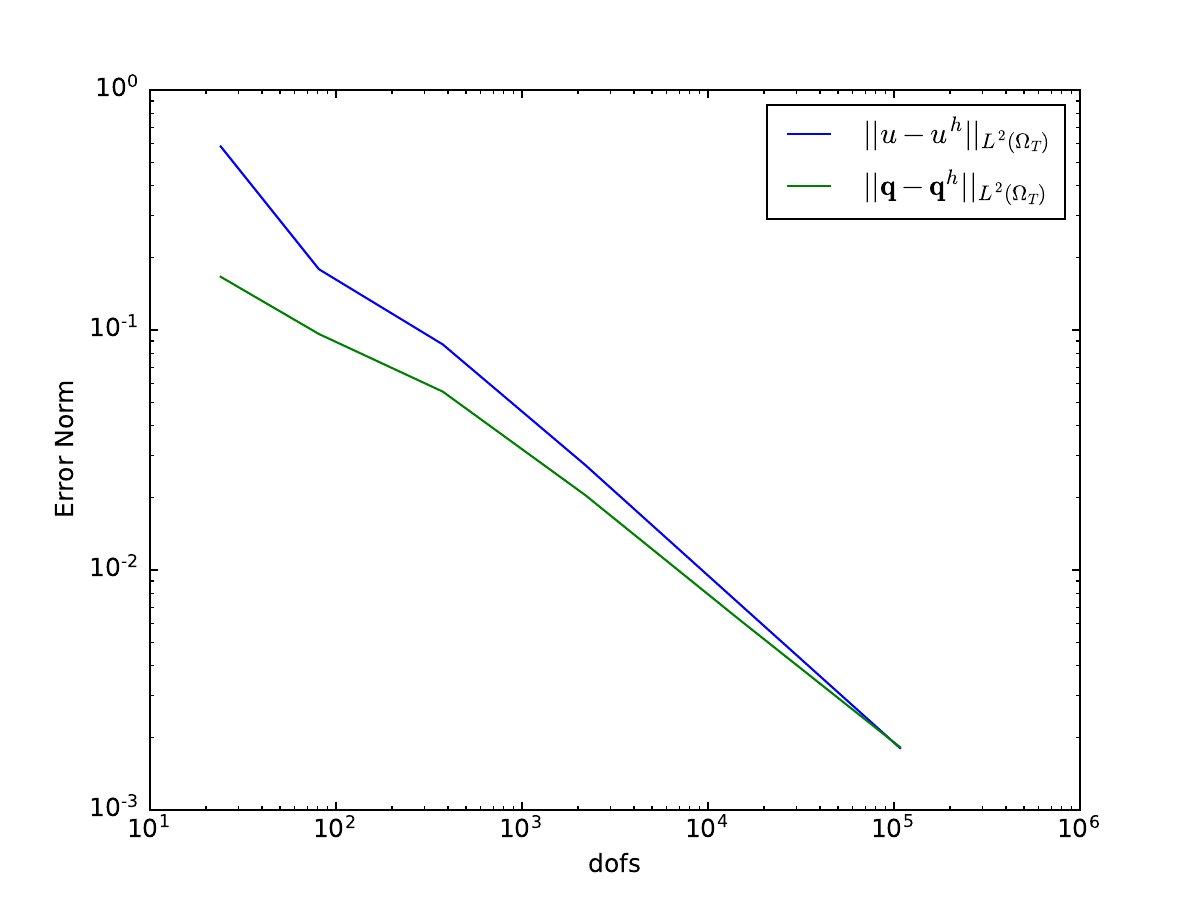}}
    \subfigure[ \label{fig:asympt_convp2_spacetime} Quadratic polynomial approximations.  ]{\centering
 \includegraphics[width=0.51\textwidth]{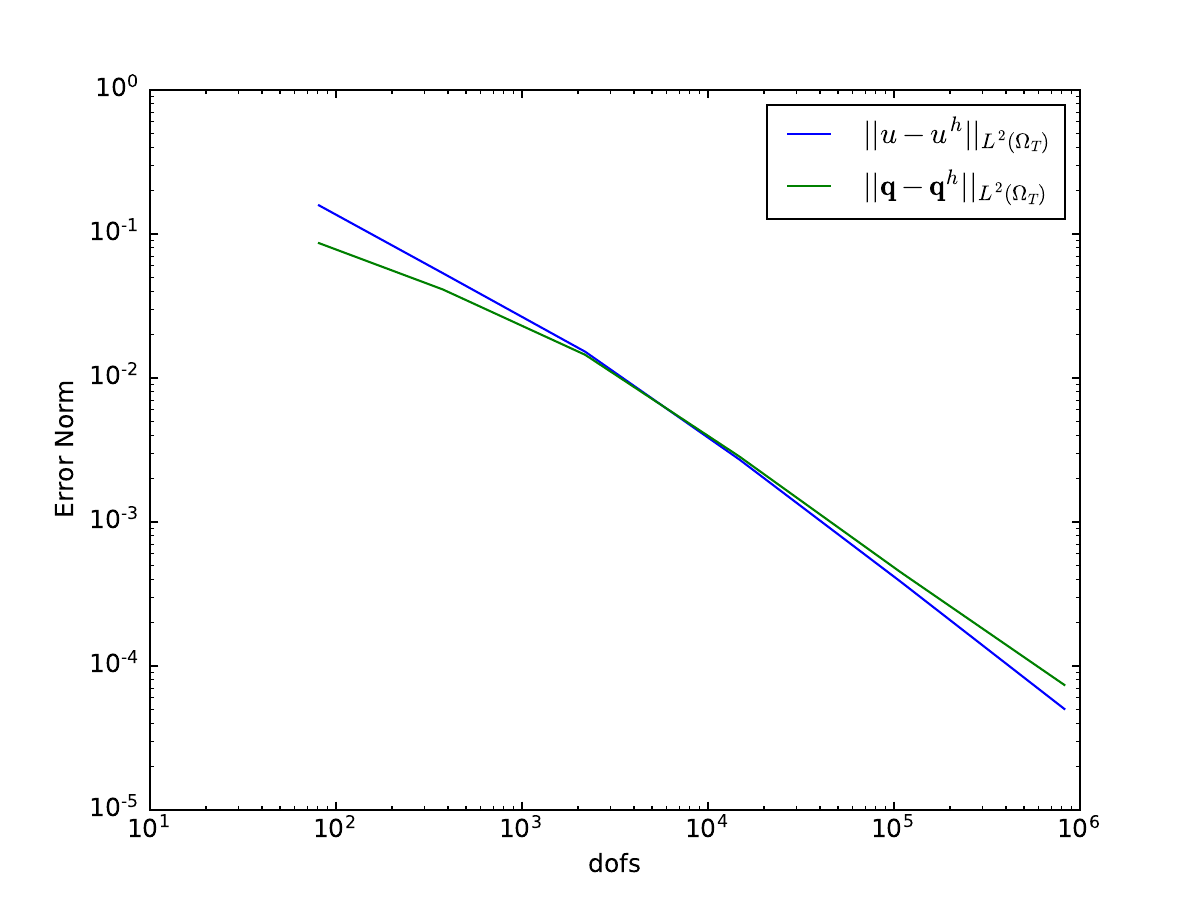}}
  \caption{\label{fig:asympt_conv_spacetime} Convergence histories for the space-time convergence study.}
\end{figure}

In Figure~\ref{fig:asympt_conv_spacetime} the convergence plots for linear and quadratic polynomial degrees for the space-time approach are shown. 
%
In Figure~\ref{fig:asympt_conv_spacetime}, we plot error norms versus the number of degrees of freedom $N$,
which increases at $\mathcal{O}(h^{-2})$, i.e., the $h-$convergence rates of the FE approximations 
can be extracted from these by a simple adjustment. For the case of $\norm{u-u^h}{\SLTOT}$, we get $\mathcal{O}(N^{-1}) = \mathcal{O}(h^{2}) =\mathcal{O}(h^{p+1})$ order of convergence. 
The observed rates for $\norm{\qq-\qq^h}{\SLTOT}$ are slightly lower, whereas the energy error converges at the expected rates of $\mathcal{O}(h^{p})$. 
In error bounds for the AVS-FE method applied to a second order PDE, see, e.g,~\cite{valseth2021stable}, it is only guaranteed that the energy norm~\eqref{eq:energy_norm} and the error in the norm on $\UUUT$ converges at 
$\mathcal{O}(h^{p})$. 

\begin{figure}[h!]
	\subfigure[ \label{fig:GAr0p1} Linear polynomial approximations, $\rho_\infty=0$. ]{\centering
		\includegraphics[width=0.5\textwidth]{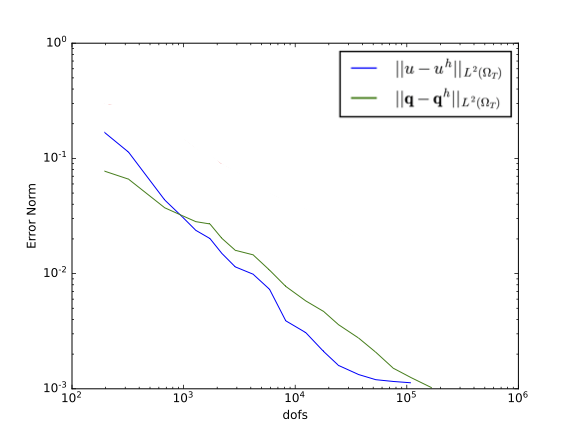}}
	\subfigure[ \label{fig:GAr9p1} Linear polynomial approximations, $\rho_\infty=0.9$.  ]{\centering
		\includegraphics[width=0.5\textwidth]{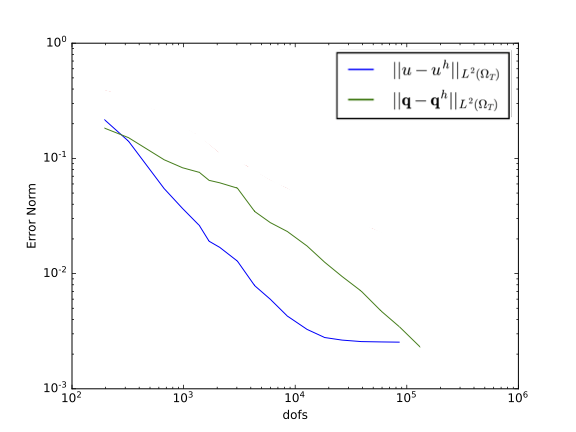}}
		\subfigure[ \label{fig:GAr0p2} Quadratic polynomial approximations, $\rho_\infty=0$.  ]{\centering
		\includegraphics[width=0.5\textwidth]{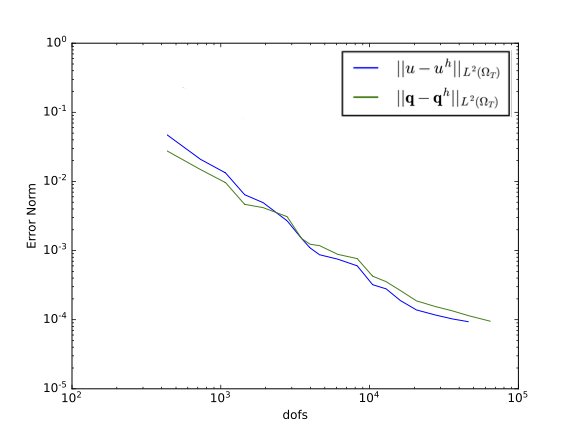}}
		\subfigure[ \label{fig:GAr9p2} Quadratic polynomial approximations, $\rho_\infty=0.9$.  ]{\centering
		\includegraphics[width=0.5\textwidth]{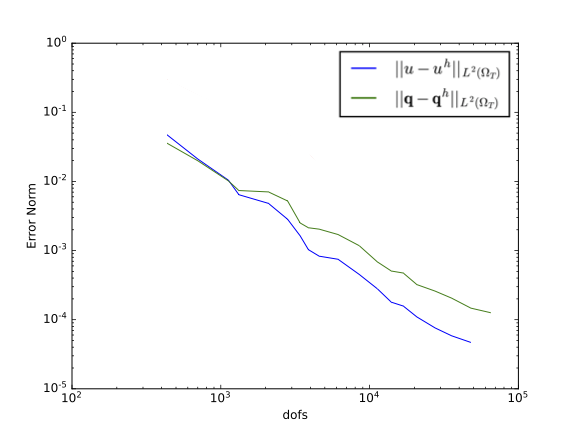}}
	\caption{\label{fig:GAerror} Convergence study of the solution obtained using the generalized-$\alpha$ method for time discretization at final time $T=0.5\,s$. }
\end{figure}

Analogously, in Figure~\ref{fig:GAerror}, the convergence plots for generalized-$\alpha$ are 
presented for to study the convergence of the method at the final time $T=0.5s$  with time step of $\tau=10^{-3}$. 
The observed rates of convergence in Figure~\ref{fig:GAerror} are the optimal rates expected from the polynomial 
approximations employed. Note that the $L^2$ errors in the base variable $u$ becomes flat near the end of the refinement process as the temporal discretization error becomes dominant.
Comparison of the results in Figures~\ref{fig:asympt_conv_spacetime} 
and~\ref{fig:GAerror} for the two methods reveal that the number of degrees of freedom 
is significantly larger for the space-time approach.

\subsubsection{$H-div$ Conforming  Basis Functions }
To complete our numerical verifications we consider the generalized-$\alpha$ system~\eqref{eq:mix_par} and use Raviart-Thomas basis functions for the flux $q^h$. Following the know results from e.g.,~\cite{BrezziMixed}, the Raviart-Thomas functions are of order $p-1$, where $p$ is the order of the approximations for $u^h$.
We also use discontinuous Raviart-Thomas bases for the vector valued error representation function of order $p$ and the scalar valued function of the same order $p$. 

\begin{figure}[h!]
	\subfigure[  Uniform mesh refinements. ]{\centering
		\includegraphics[width=0.5\textwidth]{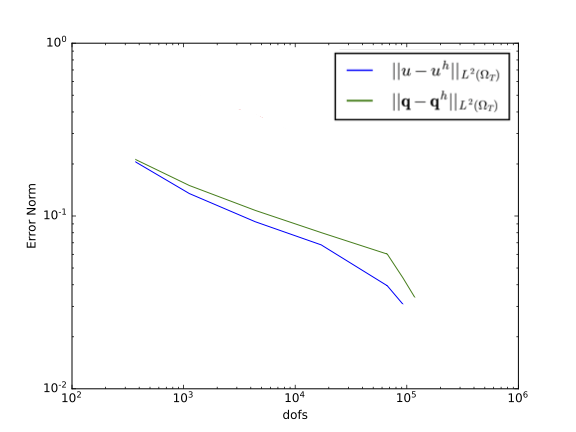}}
	\subfigure[ Adaptive mesh refinements.   ]{\centering
		\includegraphics[width=0.5\textwidth]{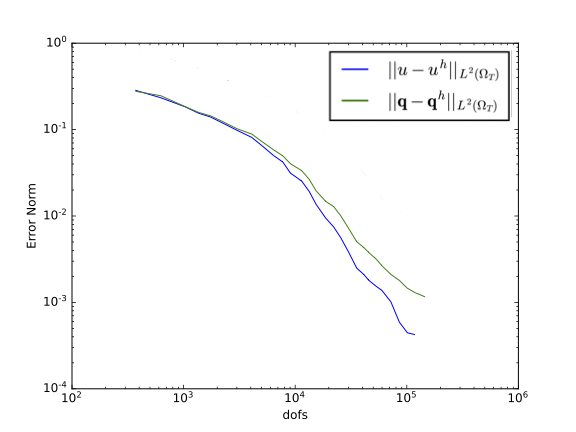}}
	\caption{\label{fig:GARTerror} Convergence study of the solution obtained using the generalized-$\alpha$ approach $T=1.0$ using fully conforming FE basis functions. }
\end{figure}

We again consider the same Eriksson-Johnson problem with $T_{final}=1.0s$, set $\epsilon=1\times 10^{-3}$, $p=2$, $\rho_\infty = 0.9$ and perform both uniform and mesh refinements. In Figure~\ref{fig:GARTerror}, we present the corresponding convergence histories. Clearly, for the strongly convection-dominated case considered, the uniform refinements are not an optimal choice. However, the adaptive refinement scheme performs significantly better and is able to reduce the considered errors approximately two orders of magnitude. 

\subsubsection{Comparison Between Space-Time and Time Stepping}
\label{sec:compare_two}

As the space-time and time-stepping methods are fundamentally different, a comparison between the two methods is not trivial. Comparison of accuracy of the two methods is not straightforward to compare, as the 
errors reported in Figure~\ref{fig:asympt_conv_spacetime} are global for the full space-time domain and the errors in Figure~\ref{fig:GAerror} are at the final time step. Furthermore, the computational cost is distributed differently in the two methods. To provide a heuristic comparison between the two methods, we compare the error at the final time $T=0.5$ for the case of $p=1$ with the problem setup from Section \ref{sec:convergence}. In the space-time approach the initial mesh consists of six uniform space-time tetrahedrons whereas in the GA method it consists of two triangular elements. In the GA method we set $\rho_\infty = 0.9$ and perform 5 time steps. We perform uniform to the initial mesh and compute the errors in the space-time approach at the final time step and plot them alongside the final time error from generalized-$\alpha$ against the (2 dimensional) element size $h$ at the final time in Figure~\ref{fig:ErrorsAdvDiff}.   It is interesting to observe that the errors in both methods shown in this figure are nearly identical. In terms of computational time, the space-time approach required 75 seconds whereas the GA method took 49 seconds. In both cases the experiments were performed on a 2022 MacBook pro with the Apple M2 chip.  

\begin{figure}[h!]
\centering
\begin{tikzpicture}
\begin{axis}[small,
width=1.45*\width,
height=\height,
ymode=log,
xmode=log,
xlabel=h,
ylabel=$\norm{u-u^h}{\SLTO}$,
grid=both,
    legend style={font=\small,at={(0.5,1.03)},anchor=south},
    legend entries={Space-time, Time stepping},
    legend columns=3
]      
\addplot[
color=red,dashdotted,line width=1pt, 
mark=*,mark size=1.5pt,    
]
table[x=h,y=L2errU]{new_figs/Eirik/Spacetime_final.txt};
\addplot[
color=blue,dotted,line width=1pt, 
mark=*,mark size=1.5pt,    
]
table[x=h,y=L2errU]{new_figs/Eirik/timestep_final.txt};
\end{axis}
\end{tikzpicture}
\caption{Convergence at the final time  $T=0.5$ for increasingly fine uniform meshes. \label{fig:ErrorsAdvDiff}}

\end{figure}

\clearpage

\subsection{Shock Problem}
\label{sec:shock}
\begin{figure}[h!]
\hspace*{-1.95in}\subfigure[ \label{fig:shock02s_mesh} Solution $u^h$ at $t=0.1s$. ]{\centering
 \includegraphics[width=0.95\textwidth]{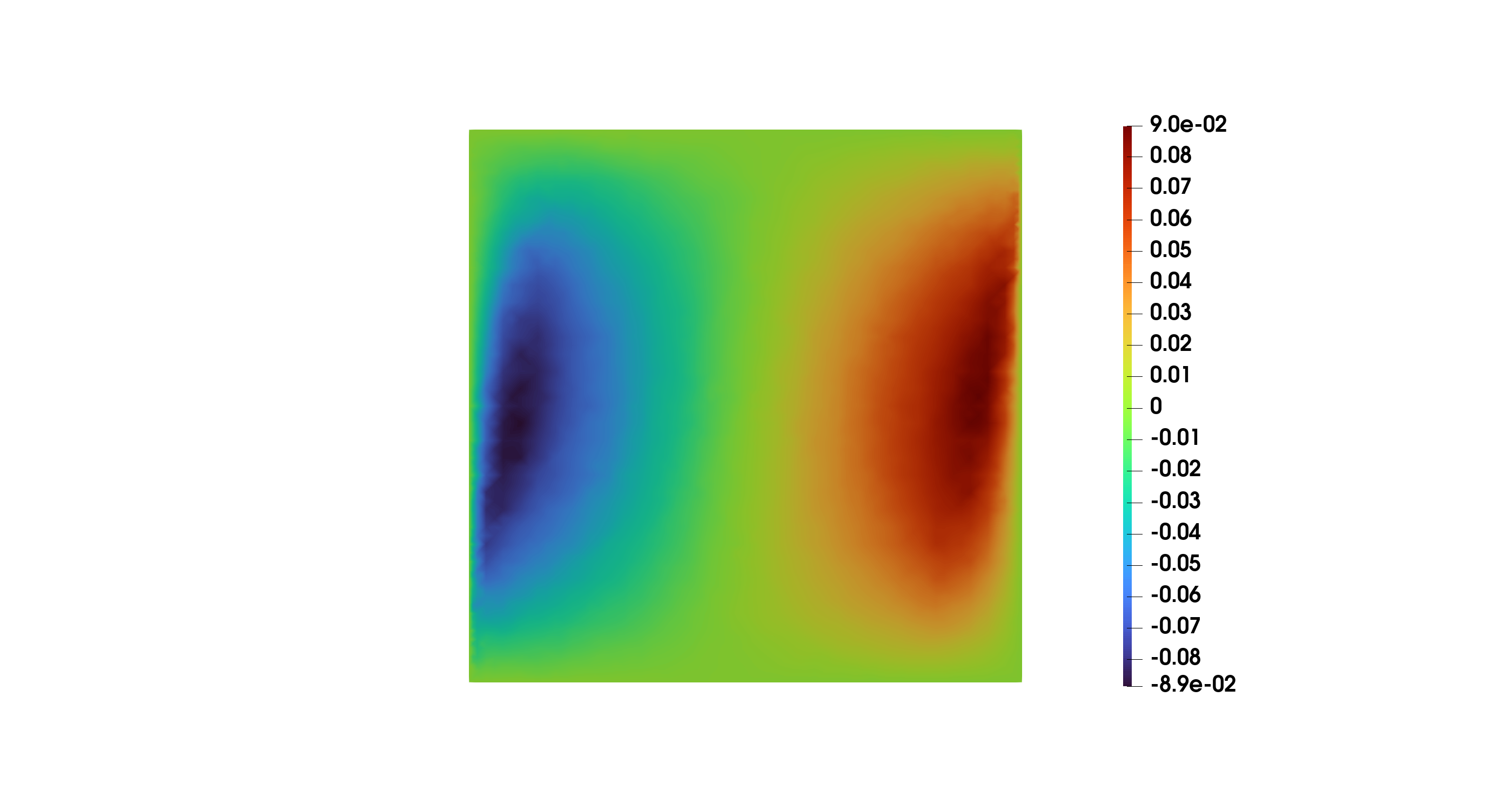}}
\hspace*{-2.85in}    \subfigure[ \label{fig:shock05s} Solution $u^h$ at $t=0.2s$ with final adapted mesh. ]{\centering
 \includegraphics[width=0.95\textwidth]{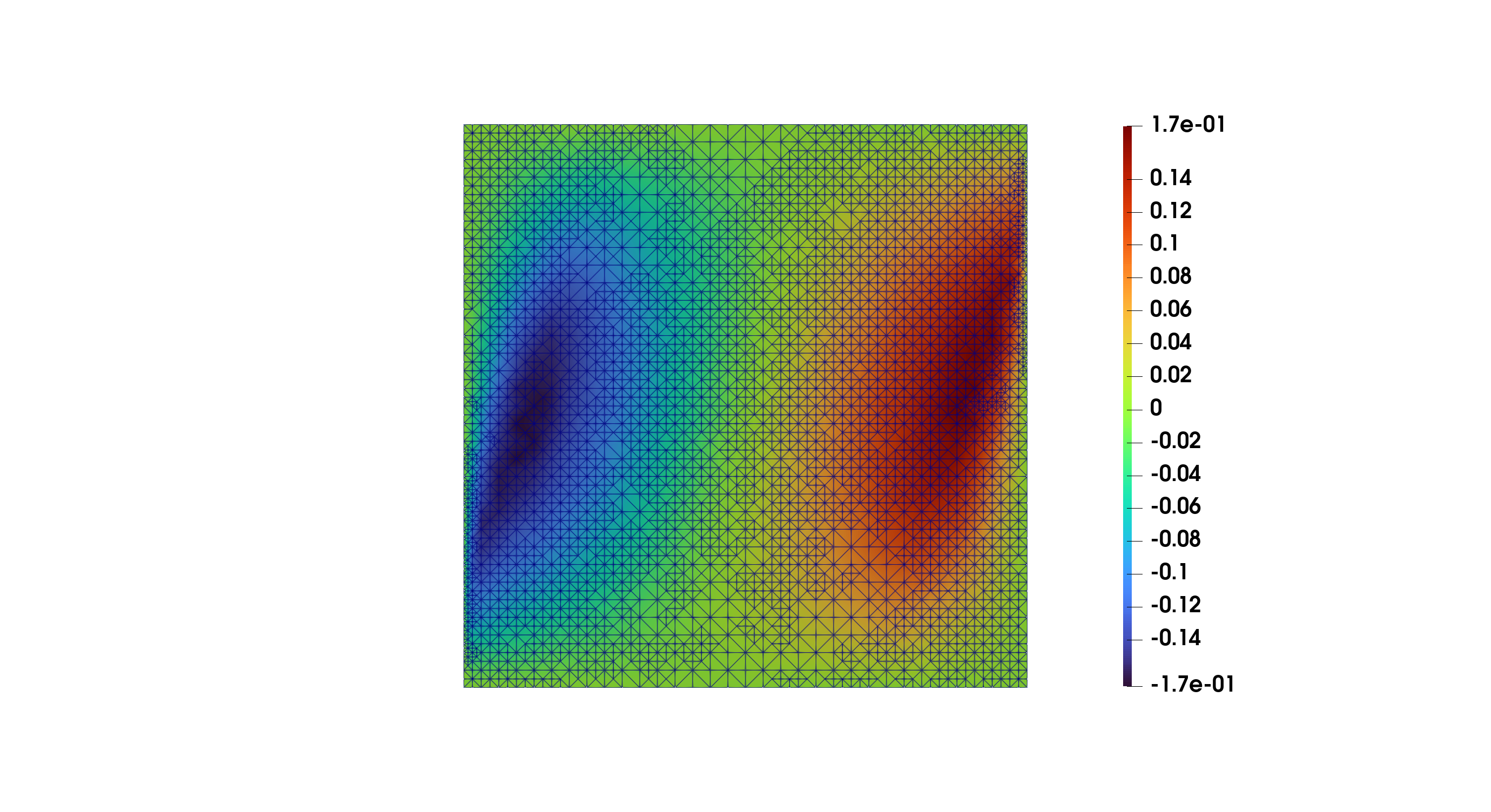}}
  \caption{\label{fig:shock_figure_1} AVS-FE approximations of the shock problem, i.e., \eqref{eq:num_results_homo_conv_diff} with parameters from~\eqref{eq:shock_params}.}
\end{figure}

\begin{figure}[h!]
\hspace*{-1.95in}\subfigure[ \label{fig:shock1s} Solution $u^h$ at $t=1.0s$. ]{\centering
 \includegraphics[width=0.96\textwidth]{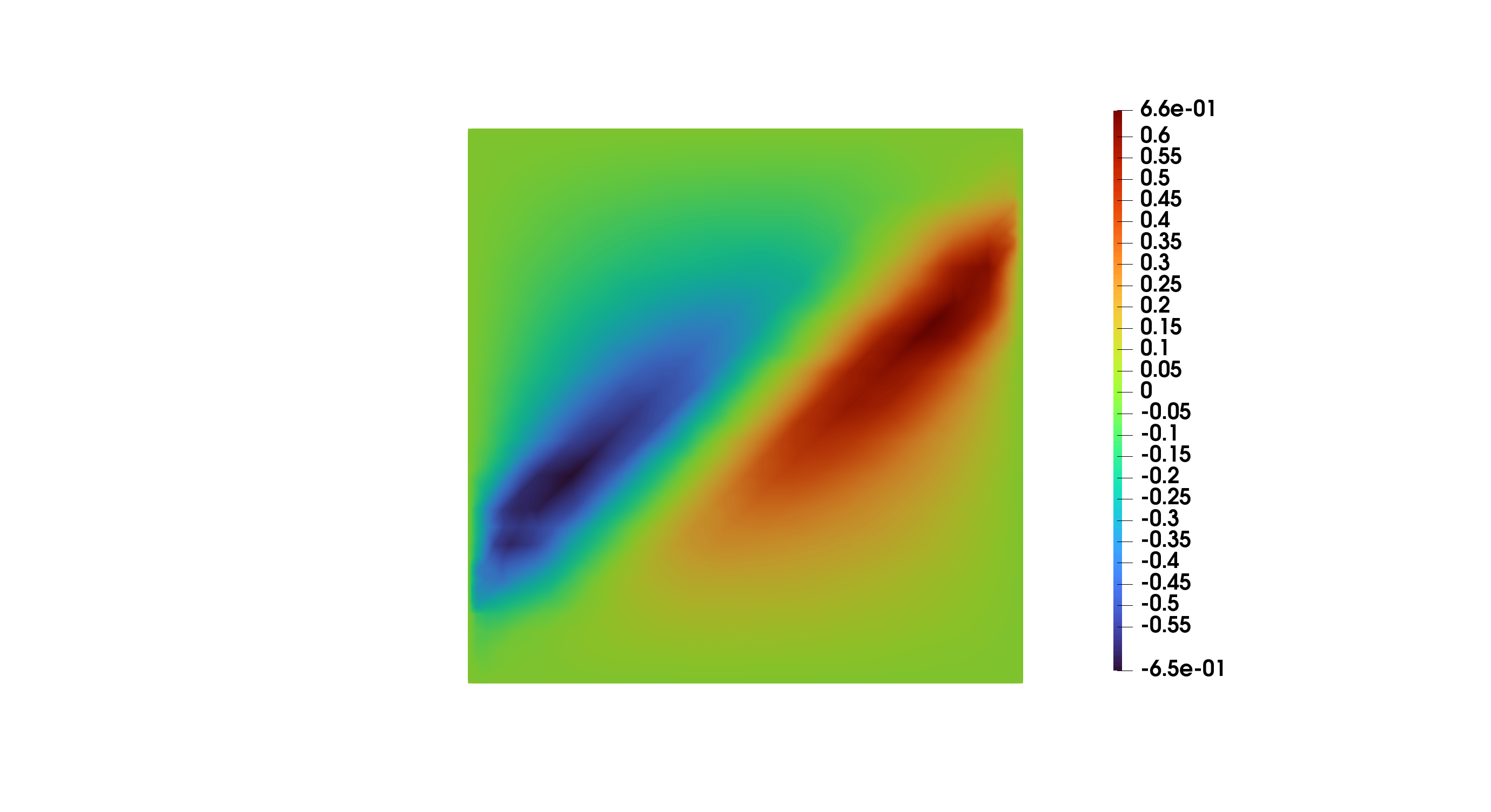}}
\hspace*{-2.85in}    \subfigure[ \label{fig:shock135s_mesh} Solution $u^h$ at $t=1.35s$ with final adapted mesh. ]{\centering
 \includegraphics[width=0.95\textwidth]{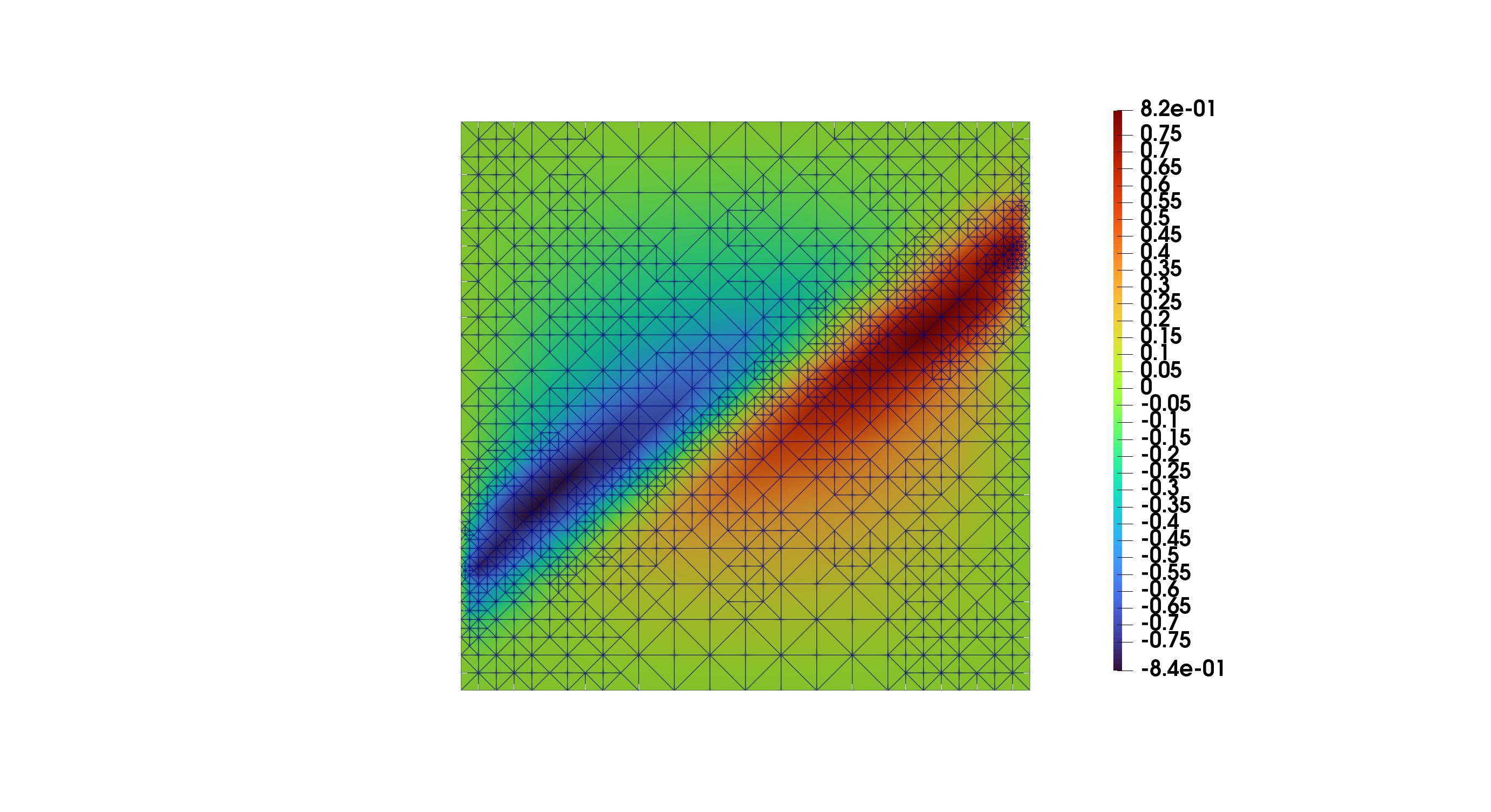}}
  \caption{\label{fig:shock_figure_2} AVS-FE approximations of the shock problem, i.e., \eqref{eq:num_results_homo_conv_diff} with parameters from~\eqref{eq:shock_params}.}
\end{figure}

\begin{figure}[h!]
\hspace*{-1.95in}\subfigure[ \label{fig:shock2s} Solution $u^h$ at $t=2.0s$. ]{\centering
 \includegraphics[width=0.995\textwidth]{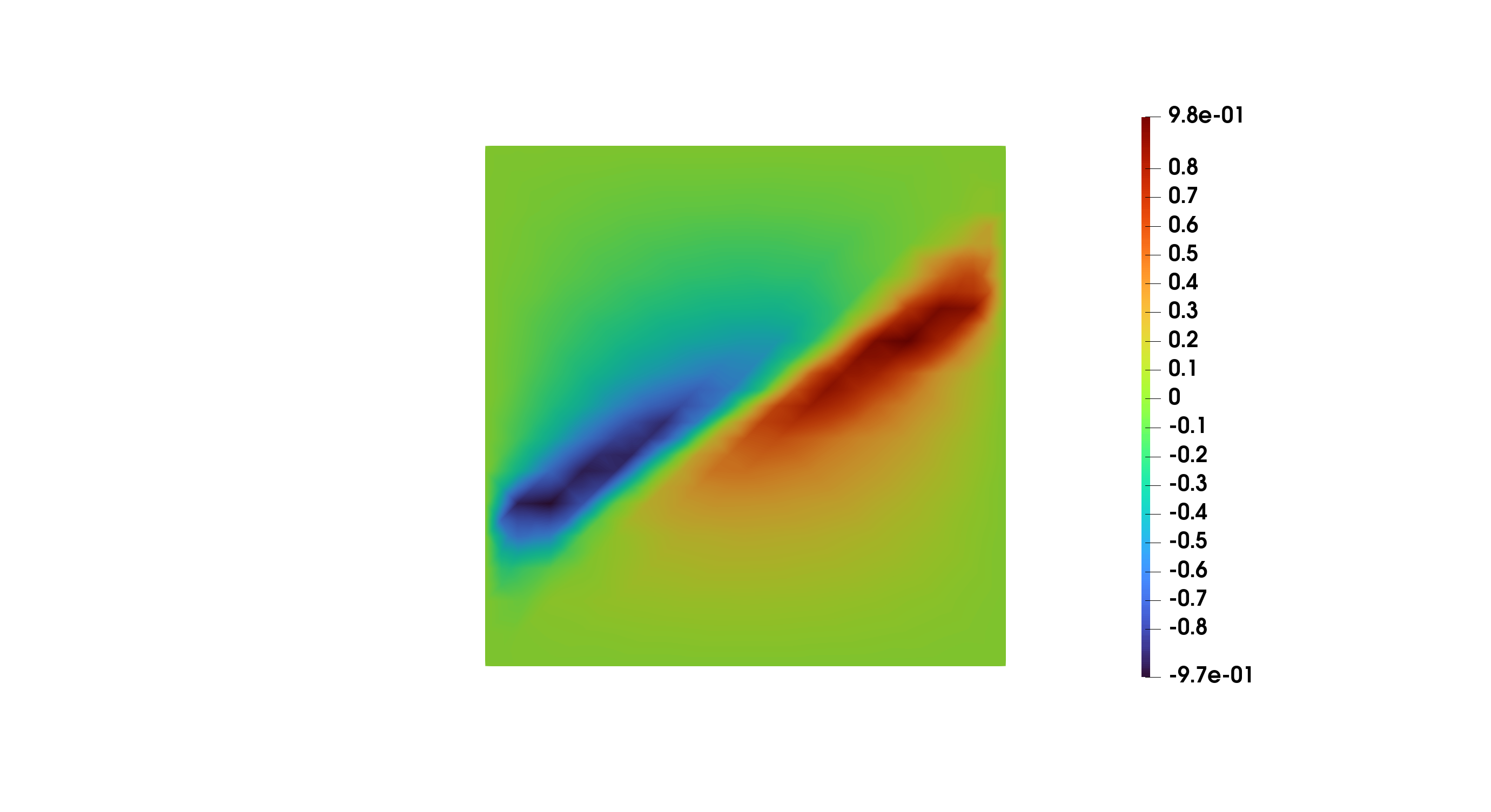}}
\hspace*{-2.85in}    \subfigure[ \label{fig:shock250s_mesh} Solution $u^h$ at $t=2.5s$ with final adapted mesh. ]{\centering
 \includegraphics[width=0.95\textwidth]{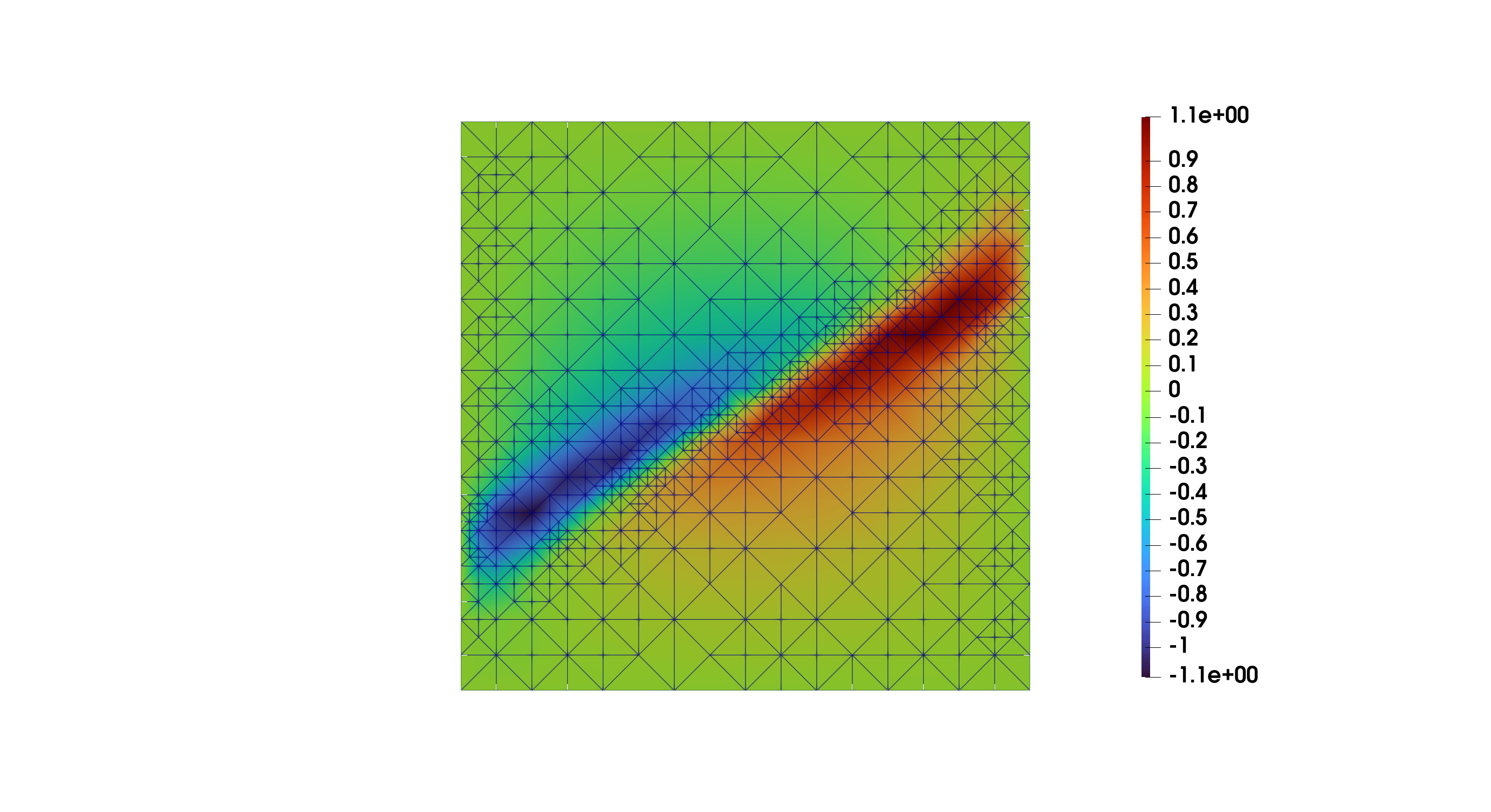}}
  \caption{\label{fig:shock_figure_3} AVS-FE approximations of the shock problem, i.e., \eqref{eq:num_results_homo_conv_diff} with parameters from~\eqref{eq:shock_params}.}
\end{figure}

As a final numerical verification, we present a consideration of~\eqref{eq:num_results_homo_conv_diff} in which the solution behaves as two shocks 
traveling through the space-time domain while rotating about the origin. Furthermore, the choices we make for the 
problem parameters are such that the interface of the shock is skewed and rotates 
in the space-time domain as $t\rightarrow T_{final}$. Thus, we have the following 
choices:
\begin{equation}\label{eq:shock_params}  
\begin{array}{rl}
\ds  \bb = & \{-x + 2y, 0 \}^T,  \\
\ds  u_{0} = & 0, \\
\ds  \epsilon = & 10^{-3}, \\ 
\ds  u_{0} = & 0, \\
\ds  f = & -2x \epsilon + x(1-y^2). \\
\ds  T_{final}= & 2.50s \\
\ds  \Omega = & (-1, 1 ) \times( -1,1 ) \\
 \end{array} 
\end{equation}
For this particular problem, we present the time slice approach in which 
we perform mesh adaptations between each slice and we apply linear polynomial approximations for the trial functions.  
Experience has shown that the slice containing the initial condition is critical to the proper resolution of the space-time 
process. Thus, we consider the case of three space-time slices, the first from $0s$ to $0.2s$ and the final two of equal size from $0.2s$ to $2.5s$. In Figures~\ref{fig:shock_figure_1},~\ref{fig:shock_figure_2},  and~\ref{fig:shock_figure_3} we present the AVS-FE solution for the 
base variable at different time steps. As expected, two shock-waves originate at the boundaries of $x=\pm 1$, and as 
time progress, the two waves approach the center of the domain while rotating.
The adaptively refined meshes shown in Figures~\ref{fig:shock05s}, ~\ref{fig:shock135s_mesh},  and~\ref{fig:shock250s_mesh} 
  (the final times of each slice) show that the mesh refinements are focused at the interfaces of  the shocks, further indicating the applicability of the built-in error indicators.
%
%
%
%
%
%
%

\section{Conclusions}
\label{sec:conclusions}
The AVS-FE method is a Petrov-Galerkin method which uses classical 
continuous FE trial basis functions, while the test space consist of functions that are 
discontinuous across element edges. 
This broken topology in the test space allows us 
to employ the DPG philosophy and introduce an equivalent saddle point problem which we implement using 
high level FE solvers.
We have introduced two distinct approaches to transient problems using the AVS-FE method. First, we 
take a space-time approach in which the entire space-time domain in discretized using finite elements,
and second, using the method of lines to discretize the spatial domain independently to deliver a semi-discretized system. Then, using a time-marching method, we obtain a fully discrete system.

The space-time method allows us to exploit the unconditional stability of the AVS-FE method and perform
a single global solve governing the FE approximation. As the AVS-FE approximations computed from the saddle point system~\eqref{eq:min_prob} come with built-in error indicators, we are capable of utilizing mesh adaptive strategies in space and time. 
 In an effort to control  the computational cost of the space-time approach
in solving the global system of equations, we consider a time slice approach. Here, the space-time domain 
is partitioned into finite sized space-time slices on which we employ the AVS-FE method. The advantage 
here is that the size of the global system is reduced and we are able to employ mesh adaptive 
strategies on each slice.

The method of lines, in which we use the AVS-FE method for 
the spatial discretization and a generalized-$\alpha$ method to derive a fully-discretized system.
In this case, the discrete stability in the temporal domain is ensured by the generalized-$\alpha$ method leading to 
highly efficient stable FE computations. We show that the AVS-FE method uses a corresponding norm as a function of the time-step. Another distinguishing feature of this method is that due to the influence of the initial data on the accuracy of the solution, we find a stable approximation for $\frac{\partial u}{\partial t}$ at the initial time. Accordingly, at each time step, one requires to solve a system with a smaller number of degrees of freedom in comparison with the space-time approach. 

Numerical verifications for several cases of the transient convection-diffusion IBVP show that 
both methods exhibit optimal asymptotic convergence behavior as well as similar norms of the numerical 
approximation error. For degrees of approximation above $2$, the space-time approach becomes 
more accurate as it is not limited to the second-order accuracy of the generalized-$\alpha$ method. 
However, we do not advocate one method over the other but we point out these differences for 
potential users as their available computational resources will likely dictate which approach to use.
For both cases, we present additional numerical verifiactions highlighting the adaptive mesh 
refinement capabilities. 
In future efforts, we expect to pursue alternative error estimators and indicators as well 
as the AVS-FE approximation of challenging transient physical phenomena.
The use of basis functions that are of higher order regularity, e.g, as in \cite{los2020isogeometric} is another potential direction of future research efforts.

\section*{Acknowledgements}
Authors Valseth and Dawson have been supported by the United States National Science Foundation - NSF PREEVENTS Track 2 Program,under  NSF Grant Number  1855047.
Authors Valseth and Romkes 
have been supported by the United States National Science Foundation - NSF CBET Program,
under  NSF Grant Number 1805550.
The authors gratefully acknowledge the assistance of Austin Kaul of the Department of Mechanical Engineering at South Dakota School of Mines and Technology in performing the numerical 
verifications of Section~\ref{sec:convergence} and~\ref{sec:shock}. Finally, the authors would also like to gratefully acknowledge the use of the “ADCIRC” and “DMS21031” allocations at the Texas Advanced Computing Center at the University of Texas at Austin.

\bibliographystyle{elsarticle-num}
 \bibliography{references_eirik}

\begin{thebibliography}{10}
\expandafter\ifx\csname url\endcsname\relax
  \def\url#1{\texttt{#1}}\fi
\expandafter\ifx\csname urlprefix\endcsname\relax\def\urlprefix{URL }\fi
\expandafter\ifx\csname href\endcsname\relax
  \def\href#1#2{#2} \def\path#1{#1}\fi

\bibitem{CaloRomkesValseth2018}
V.~M. Calo, A.~Romkes, E.~Valseth, Automatic variationally stable analysis for
  {FE} computations: an introduction, in: Barrenechea G., Mackenzie J. (eds)
  Boundary and Interior Layers, Computational and Asymptotic Methods BAIL 2018,
  Springer, 2020, pp. 19--43.

\bibitem{Demkowicz4}
L.~Demkowicz, J.~Gopalakrishnan, A class of discontinuous {P}etrov-{G}alerkin
  methods. {P}art {I}: The transport equation, Computer Methods in Applied
  Mechanics and Engineering 199~(23) (2010) 1558--1572.

\bibitem{Demkowicz2}
C.~Carstensen, L.~Demkowicz, J.~Gopalakrishnan, A posteriori error control for
  {DPG} methods, SIAM Journal on Numerical Analysis 52~(3) (2014) 1335--1353.

\bibitem{Demkowicz3}
L.~Demkowicz, J.~Gopalakrishnan, Analysis of the {DPG} method for the {P}oisson
  equation, SIAM Journal on Numerical Analysis 49~(5) (2011) 1788--1809.

\bibitem{Demkowicz5}
L.~Demkowicz, J.~Gopalakrishnan, A class of discontinuous {P}etrov-{G}alerkin
  methods. {II}. {O}ptimal test functions, Numerical Methods for Partial
  Differential Equations 27~(1) (2011) 70--105.

\bibitem{Demkowicz6}
L.~Demkowicz, J.~Gopalakrishnan, A class of discontinuous {P}etrov-{G}alerkin
  methods. {P}art {III}: Adaptivity, Applied numerical mathematics 62~(4)
  (2012) 396--427.

\bibitem{deng2019high}
P.~Behnoudfar, Q.~Deng, V.~M. Calo, High-order generalized-alpha method,
  Applications in Engineering Science 4 (2020) 100021.

\bibitem{behnoudfar2019higher}
P.~Behnoudfar, Q.~Deng, V.~M. Calo, Higher-order generalized-$\alpha$ methods
  for hyperbolic problems, Computer Methods in Applied Mechanics and
  Engineering 378 (2021) 113725.

\bibitem{chung1993time}
J.~Chung, G.~Hulbert, A time integration algorithm for structural dynamics with
  improved numerical dissipation: the generalized-$\alpha$ method, Journal of
  Applied Mechanics 60~(2) (1993) 371--375.

\bibitem{Hughes1996}
T.~J.~R. Hughes, J.~R. Stewart, A space-time formulation for multiscale
  phenomena, Journal of Computational and Applied Mathematics 74 (1996)
  217--229.

\bibitem{hughes1988space}
T.~J. Hughes, G.~M. Hulbert, Space-time finite element methods for
  elastodynamics: formulations and error estimates, Computer methods in applied
  mechanics and engineering 66~(3) (1988) 339--363.

\bibitem{aziz1989continuous}
A.~K. Aziz, P.~Monk, Continuous finite elements in space and time for the heat
  equation, Mathematics of Computation 52~(186) (1989) 255--274.

\bibitem{valseth2020Cahn}
E.~Valseth, A.~Romkes, A.~R. Kaul, A stable {FE} method for the space-time
  solution of the {Cahn-Hilliard} equation, Journal of Computational Physics
  441 (2021) 110426.

\bibitem{VALSETH2020113297}
E.~Valseth, C.~Dawson, An unconditionally stable space–time {FE} method for
  the {Korteweg}–de {Vries} equation, Computer Methods in Applied Mechanics
  and Engineering 371 (2020) 113297.
\newblock \href {https://doi.org/https://doi.org/10.1016/j.cma.2020.113297}
  {\path{doi:https://doi.org/10.1016/j.cma.2020.113297}}.

\bibitem{ellis2014space}
T.~E. Ellis, L.~Demkowicz, J.~Chan, R.~D. Moser, Space-time {DPG}: Designing a
  method for massively parallel {CFD}, {ICES} report, The Institute for
  Computational Engineering and Sciences, The University of Texas at Austin
  (2014) 14--32.

\bibitem{ellis2016robust}
T.~Ellis, J.~Chan, L.~Demkowicz, Robust {DPG} methods for transient
  convection-diffusion, in: Building bridges: connections and challenges in
  modern approaches to numerical partial differential equations, Springer,
  2016, pp. 179--203.

\bibitem{roberts2015discontinuous}
N.~V. Roberts, L.~Demkowicz, R.~Moser, A discontinuous {Petrov}--{Galerkin}
  methodology for adaptive solutions to the incompressible {Navier}--{Stokes}
  equations, Journal of Computational Physics 301 (2015) 456--483.

\bibitem{roberts2014camellia}
N.~V. Roberts, Camellia: A software framework for discontinuous
  {Petrov}--{Galerkin} methods, Computers \& Mathematics with Applications
  68~(11) (2014) 1581--1604.

\bibitem{munoz2021dpg}
J.~Mu{\~n}oz-Matute, D.~Pardo, L.~Demkowicz, A {DPG}-based time-marching scheme
  for linear hyperbolic problems, Computer Methods in Applied Mechanics and
  Engineering 373 (2021) 113539.

\bibitem{munoz2022error}
J.~Mu{\~n}oz-Matute, L.~Demkowicz, D.~Pardo, Error representation of the
  time-marching {DPG} scheme, Computer methods in applied mechanics and
  engineering 391 (2022) 114480.

\bibitem{jansen2000generalized}
K.~E. Jansen, C.~H. Whiting, G.~M. Hulbert, A generalized-$\alpha$ method for
  integrating the filtered {N}avier--{S}tokes equations with a stabilized
  finite element method, Computer Methods in Applied Mechanics and Engineering
  190~(3-4) (2000) 305--319.

\bibitem{hilber1977improved}
H.~M. Hilber, T.~J. Hughes, R.~L. Taylor, Improved numerical dissipation for
  time integration algorithms in structural dynamics, Earthquake Engineering \&
  Structural Dynamics 5~(3) (1977) 283--292.

\bibitem{wood1980alpha}
W.~Wood, M.~Bossak, O.~Zienkiewicz, An alpha modification of {N}ewmark's
  method, International Journal for Numerical Methods in Engineering 15~(10)
  (1980) 1562--1566.

\bibitem{bochevLeastSquares}
P.~B. Bochev, M.~D. Gunzburger, Least-Squares Finite Element Methods, Vol. 166,
  Springer Science \& Business Media, 2009.

\bibitem{calo2019adaptive}
V.~M. Calo, A.~Ern, I.~Muga, S.~Rojas, An adaptive stabilized conforming finite
  element method via residual minimization on dual discontinuous {Galerkin}
  norms, Computer Methods in Applied Mechanics and Engineering 363 (2020)
  112891.

\bibitem{niemi2013automatically}
A.~H. Niemi, N.~O. Collier, V.~M. Calo, Automatically stable discontinuous
  {Petrov}--{Galerkin} methods for stationary transport problems: Quasi-optimal
  test space norm, Computers \& Mathematics with Applications 66~(10) (2013)
  2096--2113.

\bibitem{nagaraj2017construction}
S.~Nagaraj, S.~Petrides, L.~F. Demkowicz, Construction of {DPG Fortin}
  operators for second order problems, Computers \& Mathematics with
  Applications 74~(8) (2017) 1964--1980.

\bibitem{demkowicz2020construction}
L.~Demkowicz, P.~Zanotti, Construction of {DPG Fortin} operators revisited,
  Computers \& Mathematics with Applications 80~(11) (2020) 2261--2271.

\bibitem{alnaes2015fenics}
M.~S. Aln{\ae}s, J.~Blechta, J.~Hake, A.~Johansson, B.~Kehlet, A.~Logg,
  C.~Richardson, J.~Ring, M.~E. Rognes, G.~N. Wells, The {FEniCS} project
  version 1.5, Archive of Numerical Software 3~(100) (2015) 9--23.

\bibitem{rathgeber2017firedrake}
F.~Rathgeber, D.~A. Ham, L.~Mitchell, M.~Lange, F.~Luporini, A.~T. McRae, G.-T.
  Bercea, G.~R. Markall, P.~H. Kelly, Firedrake: automating the finite element
  method by composing abstractions, ACM Transactions on Mathematical Software
  (TOMS) 43~(3) (2017) 24.

\bibitem{demkowicz2014overview}
L.~F. Demkowicz, J.~Gopalakrishnan, An overview of the discontinuous
  {Petrov-Galerkin} method, in: Recent developments in discontinuous {Galerkin}
  finite element methods for partial differential equations, Springer, 2014,
  pp. 149--180.

\bibitem{fuentes2017coupled}
F.~Fuentes, B.~Keith, L.~Demkowicz, P.~Le~Tallec, Coupled variational
  formulations of linear elasticity and the {DPG} methodology, Journal of
  Computational Physics 348 (2017) 715--731.

\bibitem{BrezziMixed}
F.~Brezzi, M.~Fortin, Mixed and Hybrid Finite Element Methods, Vol.~15,
  Springer-Verlag, 1991.

\bibitem{brezzi1974existence}
F.~Brezzi, On the existence, uniqueness and approximation of saddle-point
  problems arising from {Lagrangian} multipliers, Publications
  math{\'e}matiques et informatique de Rennes~(S4) (1974) 1--26.

\bibitem{behnoudfar2018variationally}
P.~Behnoudfar, V.~M. Calo, Q.~Deng, P.~D. Minev, A variationally separable
  splitting for the generalized-$\alpha $ method for parabolic equations, arXiv
  preprint arXiv:1811.09351 (2018).

\bibitem{di2011mathematical}
D.~A. Di~Pietro, A.~Ern, Mathematical aspects of discontinuous Galerkin
  methods, Vol.~69, Springer Science \& Business Media, 2011.

\bibitem{ellis2016space}
T.~E. Ellis, Space-time discontinuous {Petrov}-{Galerkin} finite elements for
  transient fluid mechanics, Ph.D. thesis (2016).

\bibitem{dorfler1996convergent}
W.~D{\"o}rfler, A convergent adaptive algorithm for {P}oisson's equation, SIAM
  Journal on Numerical Analysis 33~(3) (1996) 1106--1124.

\bibitem{amestoy2006hybrid}
P.~R. Amestoy, A.~Guermouche, J.-Y. L’Excellent, S.~Pralet, Hybrid scheduling
  for the parallel solution of linear systems, Parallel computing 32~(2) (2006)
  136--156.

\bibitem{amestoy2001fully}
P.~R. Amestoy, I.~S. Duff, J.-Y. L'Excellent, J.~Koster, A fully asynchronous
  multifrontal solver using distributed dynamic scheduling, SIAM Journal on
  Matrix Analysis and Applications 23~(1) (2001) 15--41.

\bibitem{eriksson1993adaptive}
K.~Eriksson, C.~Johnson, Adaptive streamline diffusion finite element methods
  for stationary convection-diffusion problems, mathematics of computation
  60~(201) (1993) 167--188.

\bibitem{valseth2021stable}
E.~Valseth, A.~Romkes, A.~R. Kaul, C.~Dawson, A stable mixed finite element
  method for nearly incompressible linear elastostatics, International Journal
  for Numerical Methods in Engineering 122~(17) (2021) 4709--4729.

\bibitem{los2020isogeometric}
M.~{\L}o{\'s}, J.~Munoz-Matute, I.~Muga, M.~Paszy{\'n}ski, Isogeometric
  residual minimization method {(iGRM)} with direction splitting for
  non-stationary advection--diffusion problems, Computers \& Mathematics with
  Applications 79~(2) (2020) 213--229.

\end{thebibliography}
\end{document}